\documentclass[12pt]{amsart}%
\usepackage{amsmath}
\usepackage{amsfonts}
\usepackage{amssymb}
\usepackage{amstext}
\usepackage{graphicx}%
\providecommand{\U}[1]{\protect \rule{.1in}{.1in}}
\providecommand{\U}[1]{\protect \rule{.1in}{.1in}}
\providecommand{\U}[1]{\protect \rule{.1in}{.1in}}
\newtheorem{theorem}{Theorem}[section]
\newtheorem{lemma}{Lemma}[section]
\newtheorem{proposition}{Proposition}[section]
\newtheorem{corollary}{Corollary}[section]

\newtheorem{definition}{Definition}[section]

\numberwithin{equation}{section}

\theoremstyle{remark}
\newtheorem{remark}{Remark}[section]
\numberwithin{equation}{section}
\begin{document}
\title[CR Sub-Laplacian Comparison and Liouville-type Theorem ]{CR Sub-Laplacian Comparison and Liouville-type Theorem in a Complete
Noncompact Sasakian Manifold }
\author{$^{\ast}$Shu-Cheng Chang$^{1}$}
\address{$^{1}$Department of Mathematics and Taida Institute for Mathematical Sciences
(TIMS), National Taiwan University, Taipei 10617, Taiwan\\
Current Address : Yau Mathematical Sciences Center, Tsinghua University,
Beijing, China}
\email{scchang@math.ntu.edu.tw }
\author{$^{\ast}$Ting-Jung Kuo$^{2}$}
\address{$^{2}$Department of Mathematics, National Taiwan Normal University, Taipei
11677, Taiwan }
\email{tjkuo1215@ntnu.edu.tw}
\author{$^{\ast}$Chien Lin$^{3}$}
\address{$^{3}$Department of Mathematics, National Tsing Hua University, Hsinchu 30013, Taiwan}
\email{r97221009@ntu.edu.tw}
\author{Jingzhi Tie$^{4}$}
\address{$^{4}$Department of Mathematics, University of Georgia, Athens, GA 30602-7403, U.S.A.}
\email{jtie@math.uga.edu}
\thanks{$^{\ast}$Research supported in part by the MOST of Taiwan}
\subjclass{Primary 32V05, 32V20; Secondary 53C56}
\keywords{CR Bochner formula, Subgradient estimate, Sub-Laplacian comparison theorem,
Liouvile-type theorem, Pseudohermtian Ricci, Pseudohermitian torsion, Ricatti
inequality, Sasakain manifold.}

\begin{abstract}
In this paper, we first obtain the sub-Laplacian comparison theorem in a
complete noncompact pseudohermitian manifold of vanishing torsion (i.e.
Sasakian manifold). Secondly, we derive the sub-gradient estimate for positive
pseudoharmonic functions in a complete noncompact pseudohermitian manifold
which satisfies the CR sub-Laplacian comparison property. It is served as the
CR analogue of Yau's gradient estimate. As a consequence, we have the natural
CR analogue of Liouville-type theorems in a complete noncompact Sasakian
manifold of nonnegative pseudohermitian Ricci curvature tensors.

\end{abstract}
\maketitle

\section{Introduction}

In \cite{y1} and \cite{cy}, S.-Y. Cheng and S.-T. Yau derived a well known
gradient estimate for positive harmonic functions in a complete noncompact
Riemannian manifold.

\begin{proposition}
(\cite{y1}, \cite{cy}) Let $M$ be a complete noncompact Riemannian
$m$-manifold with Ricci curvature bounded from below by $-K\ (K\geq0).$ If
$\ u\left(  x\right)  $ is a positive harmonic function on $M,$ then there
exists a positive constant $C=C(m)$ such that
\begin{equation}
|\nabla f(x)|^{2}\leq C(\sqrt{K}+\frac{1}{R})\label{112b}%
\end{equation}
on the ball $B\left(  R\right)  $ with $f(x)=\ln u(x).$ As a consequence, the
Liouville theorem holds for complete noncompact Riemannian $m$-manifolds of
nonnegative Ricci curvature.
\end{proposition}

In this paper, by modifying the arguments of \cite{y1}, \cite{cy} and
\cite{ckl}, we derive a sub-gradient estimate for positive pseudoharmonic
functions in a complete noncompact pseudohermitian $(2n+1)$-manifold
$(M,J,\theta)$ of vanishing pseudohermitian torsion (i.e. Sasakian manifold)
which is an odd dimensional counterpart of K\"{a}hler geometry. It is served
as the CR version of Yau's gradient estimate. As a consequence, we prove that
the CR analogue of Liouville-type theorem holds for complete noncompact
Sasakian manifolds of nonnegative pseudohermitian Ricci curvature.

We first define $Ric$ and $Tor$ on $T_{1,0}$ by
\begin{equation}
Ric(X,Y)=R_{\alpha \bar{\beta}}X^{\alpha}Y^{\bar{\beta}} \label{2011D}%
\end{equation}
and
\begin{equation}%
\begin{array}
[c]{c}%
Tor(X,Y)=i\sum_{\alpha,\beta}(A_{\bar{\alpha}\bar{\beta}}X^{\bar{\alpha}%
}Y^{\bar{\beta}}-A_{\alpha \beta}X^{\alpha}Y^{\beta}).
\end{array}
\label{2011C}%
\end{equation}
Here $X=X^{\alpha}Z_{\alpha}$ , $Y=Y^{\beta}Z_{\beta}$ for a frame $\left \{
T,Z_{\alpha},Z_{\bar{\alpha}}\right \}  $ of $TM\otimes \mathbb{C}$ with
$\ker \theta=\xi \mathbb{=}T_{1,0}$ $\oplus T_{0,1},$ $Z_{\alpha}\in T_{1,0}%
\ $and $Z_{\bar{\alpha}}=\overline{Z_{\alpha}}\in T_{0,1}$. $R_{\gamma}%
{}^{\gamma}{}_{\alpha \bar{\beta}}$ is the pseudohermitian curvature tensor,
$R_{\alpha \bar{\beta}}=R_{\gamma}{}^{\delta}{}_{\alpha \bar{\beta}}$ is the
pseudohermitian Ricci curvature tensor and $A_{\alpha \beta}$ \ is the torsion
tensor. We refer to section $2$ for more details about the notions of
pseudohermitian geometry.

In Yau's method for the proof of gradient estimates, \ one can estimate
$\Delta(\eta \left(  x\right)  |\nabla f(x)|^{2})$ for a nonegative cut-off
function $\eta \left(  x\right)  $ on $B\left(  2R\right)  $ via Bochner
formula and Laplacian comparison. At the end, one has gradient estimate
(\ref{112b}) by\ applying the maximum principle to $\eta \left(  x\right)
|\nabla f(x)|^{2}.$ \ 

However in order to derive the CR subgradient estimate, one of difficulties is
to deal with the following CR Bochner formula (Lemma \ref{CR Bochner}) which
involving a term $\left \langle J\nabla_{b}\varphi,\text{ }\nabla_{b}%
\varphi_{0}\right \rangle $ that has no analogue in the Riemannian case.
\[%
\begin{array}
[c]{lll}%
\Delta_{b}\left \vert \nabla_{b}\varphi \right \vert ^{2} & = & 2\left \vert
\left(  \nabla^{H}\right)  ^{2}\varphi \right \vert ^{2}+2\left \langle
\nabla_{b}\varphi,\text{ }\nabla_{b}\Delta_{b}\varphi \right \rangle \\
& + & \left(  4Ric-2\left(  n-2\right)  Tor\right)  \left(  \left(  \nabla
_{b}\varphi \right)  _{C},\text{ }\left(  \nabla_{b}\varphi \right)
_{C}\right)  +4\left \langle J\nabla_{b}\varphi,\text{ }\nabla_{b}\varphi
_{0}\right \rangle .
\end{array}
\]
Here $\left(  \nabla^{H}\right)  ^{2},$ $\Delta_{b}$, $\nabla_{b}$ are the
subhessian, sub-Laplacian and sub-gradient respectively. We \ also denote
$\varphi_{0}=T\varphi$. \ In order to overcome this difficulty, we introduce a
real-valued function $F\left(  x,\text{ }t,\text{ }R,\text{ }b\right)
:M\times \left[  0,\text{ \thinspace}1\right]  \times \left(  0,\text{ }%
\infty \right)  \times \left(  0,\text{ }\infty \right)  \rightarrow \mathbb{R} $
by adding an extra term $\ t\eta \left(  x\right)  f_{0}^{2}\left(  x\right)  $
to $\left \vert \nabla_{b}f(x)\right \vert ^{2}$ as following%

\[
F\left(  x,\text{ }t,\text{ }R,\text{ }b\right)  =t\left(  \left \vert
\nabla_{b}f(x)\right \vert ^{2}+bt\eta \left(  x\right)  f_{0}^{2}\left(
x\right)  \right)
\]
on the Carnot-Carath\'{e}odory ball $B\left(  2R\right)  $ with a constant $b$
to be determined. In section $4$, we derive the CR subgradient estimate
(\ref{112}) and (\ref{118}) by\ applying the maximum principle to $\eta \left(
x\right)  F(x,t)$ for each fixed $t\in(0,1]$ if the CR sub-Laplacian
comparison property ( \ref{d2}) holds on $(M,J,\theta)$ which is the case when
$M$ is Sasakian (Theorem \ref{t0}).

We recall that a piecewise smooth curve $\gamma:[0,1]\rightarrow M$ is said to
be horizontal if $\gamma^{\prime}(t)\in \xi$ whenever $\gamma^{\prime}(t)$
exists. The length of $\gamma$ is then defined by
\[
l(\gamma)=\int_{0}^{1}\langle \gamma^{\prime}(t),\gamma^{\prime}(t)\rangle
_{L_{\theta}}^{\frac{1}{2}}dt.
\]
Here $\left \langle \ ,\  \right \rangle _{L_{\theta}}$ is the Levi form as in
(\ref{22a}). The Carnot-Carath\'{e}odory distance between two points $p$,
$q\in M$ is
\[
d_{c}(p,q)=\inf \{l(\gamma)\,|\, \gamma \in C_{p,q}\}
\]
where $C_{p,q}$ is the set of all horizontal curves joining $p$ and $q$. We
say $M$ is complete if it is complete as a metric space. We refer to \cite{s}
for some details. By Chow connectivity theorem \cite{cho}, there always exists
a horizontal curve joining $p$ and $q$, so the distance is finite.
Furthermore, there is a minimizing geodesic joining $p$ and $q$ so that its
length is equal to the distance $d_{c}(p,q)$.

Firstly, by applying the Ricatti inequality for sub-Laplacian of
Carnot-Caratheodory distance as in Lemma \ref{la1} and Theorem \ref{t41}, we
have the following Bishop-type sub-Laplacian comparison property in a complete
noncompact pseudohermitian $(2n+1)$-manifold of vanishing pseudohermitian
torsion tensors.

\begin{theorem}
\label{t0} Let $\left(  M,\text{ }J,\text{ }\theta \right)  $ be a complete
noncompact pseudohermitian $(2n+1)$-manifold of vanishing pseudohermitian
torsion tensors with
\[
Ric\left(  Z,\text{ }Z\right)  \geq-k\left \vert Z\right \vert ^{2}%
\]
for all $Z\in T_{1,0}$ and $k$ is an nonnegative constant. Then

(i) $n=1$%
\[
\Delta_{b}r\leq \frac{1}{r}+\sqrt{k}.
\]

(ii) $n\geq2$%
\[
\Delta_{b}r\leq \frac{2^{n}}{r}+\sqrt{2^{n}}\sqrt{k}.
\]
in the sense of distributions.
\end{theorem}

In the case of nonvanishing torsion, we make the following assumption:

\begin{definition}
\label{d2} Let $\left(  M,\text{ }J,\text{ }\theta \right)  $ be a complete
noncompact pseudohermitian $(2n+1)$-manifold with%
\begin{equation}
\left(  2Ric-(n-2)Tor\right)  \left(  Z,\text{ }Z\right)  \geq-2k\left \vert
Z\right \vert ^{2} \label{117a}%
\end{equation}
for all $Z\in T_{1,0}$, and $k$ is an nonnegative constant. We say that
$\left(  M,\text{ }J,\text{ }\theta \right)  $ satisfies the CR sub-Laplacian
comparison property if there exists a positive constant $C_{0}=C_{0}(k,n)$
\ such that
\begin{equation}
\Delta_{b}r\leq C_{0}(\frac{1}{r}+\sqrt{k}) \label{2012}%
\end{equation}

\end{definition}

We now state the following general sub-gradient estimate for positive
pseudoharmonic functions $u.$

\begin{theorem}
\label{t1} Let $\left(  M,\text{ }J,\text{ }\theta \right)  $ be a complete
noncompact pseudohermitian $(2n+1)$-manifold with%
\[
\left(  2Ric-(n-2)Tor\right)  \left(  Z,\text{ }Z\right)  \geq-2k\left \vert
Z\right \vert ^{2}%
\]
for all $Z\in T_{1,0}$, and $k\geq0.$ Furthermore, we assume that $\left(
M,\text{ }J,\text{ }\theta \right)  $ satisfies the CR sub-Laplacian comparison
property (\ref{2012}). If $u\left(  x\right)  $ is a positive pseudoharmonic
function (i.e. $\Delta_{b}u=0$) with
\begin{equation}
\lbrack \Delta_{b},T]u=0 \label{111}%
\end{equation}
on $M$. Then for each constant $b>0$, there exists a positive constant
$C_{2}=C_{2}(k)$ such that
\begin{equation}
\frac{\left \vert \nabla_{b}u\right \vert ^{2}}{u^{2}}+b\frac{u_{0}^{2}}{u^{2}%
}<\frac{\left(  n+5+2bk\right)  ^{2}}{\left(  5+2bk\right)  }\left(
k+\frac{2}{b}+\frac{C_{2}}{R}\right)  \label{118}%
\end{equation}
on the ball $B\left(  R\right)  $ of a large enough radius $R$ which depends
only on $b$,$\ k$.
\end{theorem}

\begin{remark}
It is shown that (Lemma \ref{commute})
\begin{equation}
\left[  \Delta_{b},\text{ }T\right]  u=4\operatorname{Im}[i\underset
{\alpha,\beta=1}{\overset{n}{\sum}}(A_{\bar{\alpha}\bar{\beta}}u_{\beta
}),_{\alpha}]. \label{d}%
\end{equation}
If $\left(  M,\text{ }J,\text{ }\theta \right)  $ is a complete noncompact
pseudohermitian $(2n+1)$-manifold of vanishing torsion. Then
\[
\left[  \Delta_{b},\text{ }T\right]  u=0.
\]

\end{remark}

It follows \ easily from the Theorem \ref{t0} and Theorem \ref{t1} that we
have our main results on the CR Yau's gradient estimate (\ref{119}) and
Liouville-type theorem on a complete noncompact Sasakian $(2n+1)$-manifold in
this paper.

\begin{theorem}
\label{t2018} Let $\left(  M,\text{ }J,\text{ }\theta \right)  $ be a complete
noncompact pseudohermitian $(2n+1)$-manifold of vanishing pseudohermitian
torsion and
\[
Ric\left(  Z,\text{ }Z\right)  \geq-k\left \vert Z\right \vert ^{2}%
\]
for all $Z\in T_{1,0}$, and $k\geq0.$ Let $u\left(  x\right)  $ be a positive
pseudoharmonic function. Then for each constant $b>0$, there exists a positive
constant $\overline{C}_{2}=\overline{C}_{2}(k)$ such that
\begin{equation}
\frac{\left \vert \nabla_{b}u\right \vert ^{2}}{u^{2}}+b\frac{u_{0}^{2}}{u^{2}%
}<\frac{\left(  n+5+2bk\right)  ^{2}}{\left(  5+2bk\right)  }\left(
k+\frac{2}{b}+\frac{\overline{C}_{2}}{R}\right)  \label{119}%
\end{equation}
on the ball $B\left(  R\right)  $ of a large enough radius $R$ which depends
only on $b$,$\ k$.
\end{theorem}

As a consequence, let $R$ $\rightarrow \infty$ and then $b\rightarrow \infty$
with $k=0$ in (\ref{119})$,$we have the following CR Liouville-type theorem.

\begin{corollary}
\label{c2} Let $\left(  M,\text{ }J,\text{ }\theta \right)  $ be a complete
noncompact pseudohermitian $(2n+1)$-manifold of nonnegative pseudohermitian
Ricci curvature tensors and vanishing torsion. Then any positive
pseudoharmonic function is constant.
\end{corollary}

\begin{corollary}
\label{c2b} There does not exist any positive nonconstant pseudoharmonic
function in a standard Heisenberg $(2n+1)$-manifold\emph{\ }$\left(
\mathbf{H}^{n},\text{ }\mathbf{J},\text{ }\mathbf{\theta}\right)  $.
\end{corollary}

\begin{remark}
Koranyi and Stanton (\cite{ks}) proved the Liouville theorem in $\left(
\mathbf{H}^{n},\text{ }\mathbf{J},\text{ }\mathbf{\theta}\right)  $ by a
different method.
\end{remark}

In general if the positive pseudoharmonic function $u$ does not satisfy the
condition $\left[  \Delta_{b},\text{ }T\right]  u=0$, we have the following
weak sub-gradient estimate.

\begin{theorem}
\label{t3} Let $\left(  M,\text{ }J,\text{ }\theta \right)  $ be a complete
noncompact pseudohermitian $(2n+1)$-manifold with%
\[
\left(  2Ric-(n-2)Tor\right)  \left(  Z,\text{ }Z\right)  \geq-2k\left \vert
Z\right \vert ^{2}%
\]
and
\begin{equation}
\max \left \{  \left \vert A_{\alpha \beta}\right \vert ,\  \left \vert
A_{\alpha \beta,\bar{\alpha}}\right \vert \right \}  \leq k_{1} \label{112a}%
\end{equation}
for all $Z\in T_{1,0}$ and $k\geq0,\ k_{1}>0.$ Furthermore, we assume that
$\left(  M,\text{ }J,\text{ }\theta \right)  $ satisfies the CR sub-Laplacian
comparison property. If $u\left(  x\right)  $ is a positive pseudoharmonic
function on $M$. Then there exists a small constant $b_{0}=b_{0}(n,k,k_{1})>0
$ and $C_{3}=C_{4}(k,\ k_{1},\ k_{2})$ such that for any $0<b\leq b_{0}$,
\begin{equation}
\frac{\left \vert \nabla_{b}u\right \vert ^{2}}{u^{2}}+b\frac{u_{0}^{2}}{u^{2}%
}<\frac{(n+5)^{2}}{5}\left(  k+n\left(  1+b\right)  k_{1}+\frac{2}{b}%
+\frac{C_{3}}{R}\right)  \label{112}%
\end{equation}
on the ball $B\left(  R\right)  $ of a large enough radius $R$ which depends
only on $b$.
\end{theorem}

\begin{remark}
\label{r2} By comparing the Yau's gradient estimate (\ref{112b}), we need an
extra assumption \ (\ref{112a}) to obtain the CR subgradient estimate
(\ref{112}) due to the natural of sub-Laplacian in pseudohermitian geometry.
However, we do obtain an extra estimate on the derivative of pseudoharmonic
functions $u(x)$ along the missing direction of characteristic vector field
$T$.
\end{remark}

We briefly describe the methods used in our proofs. In section $2$, we first
introduce some basic materials in a pseudohermitian $(2n+1)$-manifold. Then we
are able to get the CR Bochner-type estimate and derive some key Lemmas. In
section $3$, we give a proof of sub-Laplacian comparison theorem in a complete
noncompact pseudohermitian $(2n+1)$-manifold of vanishing pseudohermitian
torsion tensors. In section $4$, let $\left(  M,\text{ }J,\text{ }%
\theta \right)  $ be a complete noncompact pseudohermitian $(2n+1)$-manifold
with the CR sub-Laplacian comparison property, we obtain subgradient estimates
for positive pseudoharmonic functions. As a consequence, the natural analogue
of Liouville-type theorem for the sub-Laplacian holds in a complete noncompact
pseudohermitian $(2n+1)$-manifold of nonnegative pseudohermitian Ricci
curvature tensor and vanishing torsion.

\textbf{Acknowledgments. }The first author would like to express his thanks to
Prof. S.-T. Yau for the inspiration, Prof. C.-S. Lin, director of Taida
Institute for Mathematical Sciences, NTU, for constant encouragement and
supports, and Prof. J.-P. Wang for his inspiration of sublaplacian comparison
geometry. The work would be not possible without their inspirations and
supports. Part of the project was done during J.~Tie's visits to Taida
Institute for Mathematical Sciences.

\section{CR Bochner-Type Estimate}

In this section, we derive some key lemmas. In particular, we obtain the CR
Bochner-type estimate as in Lemma \ref{Bochner inequality}. We first introduce
some basic materials in a pseudohermitian $(2n+1)$-manifold (see \cite{l1},
\cite{l2} for more details).

Let $(M,\xi)$ be a $(2n+1)$-dimensional, orientable, contact manifold with
contact structure $\xi$. A CR structure compatible with $\xi$ is an
endomorphism $J:\xi \rightarrow \xi$ such that $J^{2}=-1$. We also assume that
$J$ satisfies the following integrability condition: If $X$ and $Y$ are in
$\xi$, then so are $[JX,Y]+[X,JY]$ and $J([JX,Y]+[X,JY])=[JX,JY]-[X,Y]$.

Let $\left \{  T,Z_{\alpha},Z_{\bar{\alpha}}\right \}  $ be a frame of
$TM\otimes \mathbb{C}$, where $Z_{\alpha}$ is any local frame of $T_{1,0}%
,\ Z_{\bar{\alpha}}=\overline{Z_{\alpha}}\in T_{0,1}$ and $T$ is the
characteristic vector field. Then $\left \{  \theta,\theta^{\alpha}%
,\theta^{\bar{\alpha}}\right \}  $, which is the coframe dual to $\left \{
T,Z_{\alpha},Z_{\bar{\alpha}}\right \}  $, satisfies
\begin{equation}
d\theta=ih_{\alpha \overline{\beta}}\theta^{\alpha}\wedge \theta^{\overline
{\beta}} \label{22}%
\end{equation}
for some positive definite hermitian matrix of functions $(h_{\alpha \bar
{\beta}})$. Actually we can always choose $Z_{\alpha}$ such that
$h_{\alpha \bar{\beta}}=\delta_{\alpha \beta}$; hence, throughout this note, we
assume $h_{\alpha \bar{\beta}}=\delta_{\alpha \beta}$.

The Levi form $\left \langle \ ,\  \right \rangle _{L_{\theta}}$ is the Hermitian
form on $T_{1,0}$ defined by%
\begin{equation}
\left \langle Z,W\right \rangle _{L_{\theta}}=-i\left \langle d\theta
,Z\wedge \overline{W}\right \rangle . \label{22a}%
\end{equation}
We can extend $\left \langle \ ,\  \right \rangle _{L_{\theta}}$ to $T_{0,1}$ by
defining $\left \langle \overline{Z},\overline{W}\right \rangle _{L_{\theta}%
}=\overline{\left \langle Z,W\right \rangle }_{L_{\theta}}$ for all $Z,W\in
T_{1,0}$. The Levi form induces naturally a Hermitian form on the dual bundle
of $T_{1,0}$, denoted by $\left \langle \ ,\  \right \rangle _{L_{\theta}^{\ast}%
}$, and hence on all the induced tensor bundles. Integrating the Hermitian
form (when acting on sections) over $M$ with respect to the volume form
$d\mu=\theta \wedge(d\theta)^{n}$, we get an inner product on the space of
sections of each tensor bundle. We denote the inner product by the notation
$\left \langle \ ,\  \right \rangle $. For example
\[
\left \langle u,v\right \rangle =\int_{M}u\overline{v}\ d\mu,
\]
for functions $u$ and $v$.

The pseudohermitian connection of $(J,\theta)$ is the connection $\nabla$ on
$TM\otimes \mathbb{C}$ (and extended to tensors) given in terms of a local
frame $Z_{\alpha}\in T_{1,0}$ by%

\[
\nabla Z_{\alpha}=\theta_{\alpha}{}^{\beta}\otimes Z_{\beta},\quad \nabla
Z_{\bar{\alpha}}=\theta_{\bar{\alpha}}{}^{\bar{\beta}}\otimes Z_{\bar{\beta}%
},\quad \nabla T=0,
\]
where $\theta_{\alpha}{}^{\beta}$ are the $1$-forms uniquely determined by the
following equations:%

\begin{equation}%
\begin{split}
d\theta^{\beta}  &  =\theta^{\alpha}\wedge \theta_{\alpha}{}^{\beta}%
+\theta \wedge \tau^{\beta},\\
0  &  =\tau_{\alpha}\wedge \theta^{\alpha},\\
0  &  =\theta_{\alpha}{}^{\beta}+\theta_{\bar{\beta}}{}^{\bar{\alpha}},
\end{split}
\label{23}%
\end{equation}
We can write (by Cartan lemma) $\tau_{\alpha}=A_{\alpha \gamma}\theta^{\gamma}$
with $A_{\alpha \gamma}=A_{\gamma \alpha}$. The curvature of Tanaka-Webster
connection, expressed in terms of the coframe $\{ \theta=\theta^{0}%
,\theta^{\alpha},\theta^{\bar{\alpha}}\}$, is
\[%
\begin{split}
\Pi_{\beta}{}^{\alpha}  &  =\overline{\Pi_{\bar{\beta}}{}^{\bar{\alpha}}%
}=d\omega_{\beta}{}^{\alpha}-\omega_{\beta}{}^{\gamma}\wedge \omega_{\gamma}%
{}^{\alpha},\\
\Pi_{0}{}^{\alpha}  &  =\Pi_{\alpha}{}^{0}=\Pi_{0}{}^{\bar{\beta}}=\Pi
_{\bar{\beta}}{}^{0}=\Pi_{0}{}^{0}=0.
\end{split}
\]
Webster showed that $\Pi_{\beta}{}^{\alpha}$ can be written
\begin{equation}
\Pi_{\beta}{}^{\alpha}=R_{\beta}{}^{\alpha}{}_{\rho \bar{\sigma}}\theta^{\rho
}\wedge \theta^{\bar{\sigma}}+W_{\beta}{}^{\alpha}{}_{\rho}\theta^{\rho}%
\wedge \theta-W^{\alpha}{}_{\beta \bar{\rho}}\theta^{\bar{\rho}}\wedge
\theta+i\theta_{\beta}\wedge \tau^{\alpha}-i\tau_{\beta}\wedge \theta^{\alpha}
\label{24}%
\end{equation}
where the coefficients satisfy
\[
R_{\beta \bar{\alpha}\rho \bar{\sigma}}=\overline{R_{\alpha \bar{\beta}\sigma
\bar{\rho}}}=R_{\bar{\alpha}\beta \bar{\sigma}\rho}=R_{\rho \bar{\alpha}%
\beta \bar{\sigma}},\  \  \ W_{\beta \bar{\alpha}\gamma}=W_{\gamma \bar{\alpha
}\beta}.
\]

We will denote components of covariant derivatives with indices preceded by a
comma; thus write $A_{\alpha \beta,\gamma}$. The indices $\{0,\alpha
,\bar{\alpha}\}$ indicate derivatives with respect to $\{T,Z_{\alpha}%
,Z_{\bar{\alpha}}\}$. For derivatives of a scalar function, we will often omit
the comma, for instance, $u_{\alpha}=Z_{\alpha}u,\ u_{\alpha \bar{\beta}%
}=Z_{\bar{\beta}}Z_{\alpha}u-\omega_{\alpha}{}^{\gamma}(Z_{\bar{\beta}%
})Z_{\gamma}u.$

For a real function $u$, the subgradient $\nabla_{b}$ is defined by
$\nabla_{b}u\in \xi$ and $\left \langle Z,\nabla_{b}u\right \rangle _{L_{\theta}%
}=du(Z)$ for all vector fields $Z$ tangent to contact plane. Locally
$\nabla_{b}u=\sum_{\alpha}u_{\bar{\alpha}}Z_{\alpha}+u_{\alpha}Z_{\bar{\alpha
}}$. We can use the connection to define the subhessian as the complex linear
map
\[
(\nabla^{H})^{2}u:T_{1,0}\oplus T_{0,1}\rightarrow T_{1,0}\oplus T_{0,1}%
\]
by
\[
(\nabla^{H})^{2}u(Z)=\nabla_{Z}\nabla_{b}u.\
\]
In particular,%

\[
|\nabla_{b}u|^{2}=2u_{\alpha}u_{\overline{\alpha}},\quad|\nabla_{b}^{2}%
u|^{2}=2(u_{\alpha \beta}u_{\overline{\alpha}\overline{\beta}}+u_{\alpha
\overline{\beta}}u_{\overline{\alpha}\beta}).
\]

Also
\[%
\begin{array}
[c]{c}%
\Delta_{b}u=Tr\left(  (\nabla^{H})^{2}u\right)  =\sum_{\alpha}(u_{\alpha
\bar{\alpha}}+u_{\bar{\alpha}\alpha}).
\end{array}
\]

Next we recall the following commutation relations (\cite{l1}). \ Let
$\varphi$ be a scalar function and $\sigma=\sigma_{\alpha}\theta^{\alpha}$ be
a $\left(  1,0\right)  $ form, then we have%
\begin{equation}%
\begin{array}
[c]{ccl}%
\varphi_{\alpha \beta} & = & \varphi_{\beta \alpha},\\
\varphi_{\alpha \bar{\beta}}-\varphi_{\bar{\beta}\alpha} & = & ih_{\alpha
\overline{\beta}}\varphi_{0},\\
\varphi_{0\alpha}-\varphi_{\alpha0} & = & A_{\alpha \beta}\varphi_{\bar{\beta}%
},\\
\sigma_{\alpha,0\beta}-\sigma_{\alpha,\beta0} & = & \sigma_{\alpha,\bar
{\gamma}}A_{\gamma \beta}-\sigma_{\gamma}A_{\alpha \beta,\bar{\gamma}},\\
\sigma_{\alpha,0\bar{\beta}}-\sigma_{\alpha,\bar{\beta}0} & = & \sigma
_{\alpha,\gamma}A_{\bar{\gamma}\bar{\beta}}+\sigma_{\gamma}A_{\bar{\gamma}%
\bar{\beta},\alpha},
\end{array}
\label{2010a}%
\end{equation}
and
\begin{equation}%
\begin{array}
[c]{ccl}%
\sigma_{\alpha,\beta \gamma}-\sigma_{\alpha,\gamma \beta} & = & iA_{\alpha
\gamma}\sigma_{\beta}-iA_{\alpha \beta}\sigma_{\gamma},\\
\sigma_{\alpha,\bar{\beta}\bar{\gamma}}-\sigma_{\alpha,\bar{\gamma}\bar{\beta
}} & = & ih_{\alpha \overline{\beta}}A_{\bar{\gamma}\bar{\rho}}\sigma_{\rho
}-ih_{\alpha \overline{\gamma}}A_{\bar{\beta}\bar{\rho}}\sigma_{\rho},\\
\sigma_{\alpha,\beta \bar{\gamma}}-\sigma_{\alpha,\bar{\gamma}\beta} & = &
ih_{\beta \overline{\gamma}}\sigma_{\alpha,0}+R_{\alpha \bar{\rho}}{}{}%
_{\beta \bar{\gamma}}\sigma_{\rho}.
\end{array}
\label{2010c}%
\end{equation}

Now we recall a lemma from A. Greenleaf (\cite{gr}) and also (\cite{cc2}).

\begin{lemma}
\label{CR Bochner} For a real function $\varphi$,
\begin{equation}%
\begin{array}
[c]{lll}%
\Delta_{b}\left \vert \nabla_{b}\varphi \right \vert ^{2} & = & 2\left \vert
\left(  \nabla^{H}\right)  ^{2}\varphi \right \vert ^{2}+2\left \langle
\nabla_{b}\varphi,\text{ }\nabla_{b}\Delta_{b}\varphi \right \rangle \\
& + & \left(  4Ric-2\left(  n-2\right)  Tor\right)  \left(  \left(  \nabla
_{b}\varphi \right)  _{C},\text{ }\left(  \nabla_{b}\varphi \right)
_{C}\right)  +4\left \langle J\nabla_{b}\varphi,\text{ }\nabla_{b}\varphi
_{0}\right \rangle ,
\end{array}
\label{2010}%
\end{equation}
where $\left(  \nabla_{b}\varphi \right)  _{C}=\varphi_{\bar{\alpha}}Z_{\alpha
}$ is the corresponding complex $\left(  1,\text{ }0\right)  $-vector of
$\nabla_{b}\varphi$.
\end{lemma}

\begin{lemma}
\label{Bochner inequality} For a smooth real-valued function $\varphi$ and any
$\nu>0$, we have%
\[%
\begin{array}
[c]{lll}%
\Delta_{b}\left \vert \nabla_{b}\varphi \right \vert ^{2} & \geq & 4\left(
\underset{\alpha,\beta=1}{\overset{n}{\sum}}\left \vert \varphi_{a\beta
}\right \vert ^{2}+\underset{\alpha,\beta=1,\alpha \not =\beta}{\overset{n}%
{\sum}}\left \vert \varphi_{a\bar{\beta}}\right \vert ^{2}\right)  +\frac{1}%
{n}\left(  \Delta_{b}\varphi \right)  ^{2}+n\varphi_{0}^{2}+2\left \langle
\nabla_{b}\varphi,\text{ }\nabla_{b}\Delta_{b}\varphi \right \rangle \\
& + & \left(  4Ric-2\left(  n-2\right)  Tor-\frac{4}{\nu}\right)  \left(
\left(  \nabla_{b}\varphi \right)  _{C},\text{ }\left(  \nabla_{b}%
\varphi \right)  _{C}\right)  -2\nu \left \vert \nabla_{b}\varphi_{0}\right \vert
^{2}.
\end{array}
\]

\end{lemma}%

%TCIMACRO{\TeXButton{Proof}{\proof} }%
%BeginExpansion
\proof
%EndExpansion
Since
\[%
\begin{array}
[c]{ccl}%
|(\nabla^{H})^{2}\varphi|^{2} & = & 2\sum_{\alpha,\beta=1}^{n}(\varphi
_{\alpha \beta}\varphi_{\overline{\alpha}\overline{\beta}}+\varphi
_{\alpha \overline{\beta}}\varphi_{\overline{\alpha}\beta})\\
& = & 2\sum_{\alpha,\beta=1}^{n}(|\varphi_{\alpha \beta}|^{2}+|\varphi
_{\alpha \overline{\beta}}|^{2})\\
& = & 2(\sum_{\alpha,\beta=1}^{n}|\varphi_{\alpha \beta}|^{2}+\sum
_{\substack{\alpha,\beta=1 \\ \alpha \neq \beta}}^{n}|\varphi_{\alpha
\overline{\beta}}|^{2}+\sum_{\alpha=1}^{n}|\varphi_{\alpha \overline{\alpha}%
}|^{2})
\end{array}
\]
and from the commutation relation (\ref{2010a})%
\[%
\begin{array}
[c]{ccl}%
\sum_{\alpha=1}^{n}|\varphi_{\alpha \overline{\alpha}}|^{2} & = & \frac{1}%
{4}\sum_{\alpha=1}^{n}\left(  |\varphi_{\alpha \overline{\alpha}}%
+\varphi_{\overline{\alpha}\alpha}|^{2}+\varphi_{0}^{2}\right) \\
& = & \frac{1}{4}\sum_{\alpha=1}^{n}|\varphi_{\alpha \overline{\alpha}}%
+\varphi_{\overline{\alpha}\alpha}|^{2}+\frac{n}{4}\varphi_{0}^{2}.
\end{array}
\]
It follows that%
\[%
\begin{array}
[c]{lll}%
|(\nabla^{H})^{2}\varphi|^{2} & = & 2(\sum_{\alpha,\beta=1}^{n}|\varphi
_{\alpha \beta}|^{2}+\sum_{\substack{\alpha,\beta=1 \\ \alpha \neq \beta}%
}^{n}|\varphi_{\alpha \overline{\beta}}|^{2})+\frac{1}{2}\sum_{\alpha=1}%
^{n}|\varphi_{\alpha \overline{\alpha}}+\varphi_{\overline{\alpha}\alpha}%
|^{2}+\frac{n}{2}\varphi_{0}^{2}\\
& \leq & 2(\sum_{\alpha,\beta=1}^{n}|\varphi_{\alpha \beta}|^{2}+\sum
_{\substack{\alpha,\beta=1 \\ \alpha \neq \beta}}^{n}|\varphi_{\alpha
\overline{\beta}}|^{2})+\frac{1}{2n}\left(  \Delta_{b}\varphi \right)
^{2}+\frac{n}{2}\varphi_{0}^{2}.
\end{array}
\]
On the other hand, for all $\nu>0$%
\[%
\begin{array}
[c]{lll}%
4\left \langle J\nabla_{b}\varphi,\text{ }\nabla_{b}\varphi_{0}\right \rangle  &
\geq & -4\left \vert \nabla_{b}\varphi \right \vert \left \vert \nabla_{b}%
\varphi_{0}\right \vert \\
& \geq & -\frac{2}{\nu}\left \vert \nabla_{b}\varphi \right \vert ^{2}%
-2\nu \left \vert \nabla_{b}\varphi_{0}\right \vert ^{2}.
\end{array}
\]
Then the result follows easily from Lemma \ref{CR Bochner}.
%TCIMACRO{\TeXButton{End Proof}{\endproof}}%
%BeginExpansion
\endproof
%EndExpansion

\begin{definition}
(\cite{gl}) Let $(M,J,\theta)$ be a pseudohermitian $\left(  2n+1\right)
$-manifold. We define the purely holomorphic second-order operator $Q$ by
\[
Q\varphi=2i\underset{\alpha,\beta=1}{\overset{n}{\sum}}(A_{\bar{\alpha}%
\bar{\beta}}\varphi_{\beta}),_{\alpha}.
\]

\end{definition}

By apply the commutation relations (\ref{2010a}), one obtains

\begin{lemma}
\label{commute} (\cite{gl}, \cite{ckl}) Let $\varphi \left(  x\right)  $ be a
smooth function defined on $M$. Then
\[
\Delta_{b}\varphi_{0}=\left(  \Delta_{b}\varphi \right)  _{0}+2\underset
{\alpha,\beta=1}{\overset{n}{\sum}}\left[  \left(  A_{\alpha \beta}%
\varphi_{\bar{\beta}}\right)  _{\bar{\alpha}}+\left(  A_{\bar{\alpha}%
\bar{\beta}}\varphi_{\beta}\right)  _{\alpha}\right]  .
\]
That is
\[
2\operatorname{Im}Q\varphi=\left[  \Delta_{b},\text{ }T\right]  \varphi.
\]

\end{lemma}%

%TCIMACRO{\TeXButton{Proof}{\proof} }%
%BeginExpansion
\proof
%EndExpansion
By direct computation and the commutation relation (\ref{2010a}), we have%
\[%
\begin{array}
[c]{lll}%
\Delta_{b}\varphi_{0} & = & \varphi_{0\alpha \overline{\alpha}}+\varphi
_{0\overline{\alpha}\alpha}\\
& = & \left(  \varphi_{\alpha0}+A_{\alpha \beta}\varphi_{\overline{\beta}%
}\right)  _{\overline{\alpha}}+\mathrm{conjugate}\\
& = & \varphi_{\alpha0\overline{\alpha}}+\left(  A_{\alpha \beta}%
\varphi_{\overline{\beta}}\right)  _{\overline{\alpha}}+\mathrm{conjugate}\\
& = & \varphi_{\alpha \overline{\alpha}0}+\varphi_{\overline{\alpha}\alpha
0}+2\left[  \left(  A_{\alpha \beta}\varphi_{\overline{\beta}}\right)
_{\overline{\alpha}}+\left(  A_{\bar{\alpha}\bar{\beta}}\varphi_{\beta
}\right)  _{\alpha}\right] \\
& = & \left(  \Delta_{b}\varphi \right)  _{0}+2\left[  \left(  A_{\alpha \beta
}\varphi_{\overline{\beta}}\right)  _{\overline{\alpha}}+\left(
A_{\bar{\alpha}\bar{\beta}}\varphi_{\beta}\right)  _{\alpha}\right]  .
\end{array}
\]

This completes the proof.
%TCIMACRO{\TeXButton{End Proof}{\endproof}}%
%BeginExpansion
\endproof
%EndExpansion

Let $u$ be a positive pseudoharmonic function and $f\left(  x\right)  =\ln
u\left(  x\right)  .$ Then
\[
\Delta_{b}f=-\left \vert \nabla_{b}f\right \vert ^{2}.
\]

We first define
\[
V\left(  \varphi \right)  =\underset{\alpha,\beta=1}{\overset{n}{\sum}}\left[
\left(  A_{\alpha \beta}\varphi_{\bar{\beta}}\right)  _{\bar{\alpha}}+\left(
A_{\bar{\alpha}\bar{\beta}}\varphi_{\beta}\right)  _{\alpha}+A_{\alpha \beta
}\varphi_{\bar{\beta}}\varphi_{\bar{\alpha}}+A_{\bar{\alpha}\bar{\beta}%
}\varphi_{\beta}\varphi_{\alpha}\right]  .
\]

\begin{lemma}
\label{commute f} Let $u$ be a positive pseudoharmonic function with $f=\ln u
$. Then%
\[
\Delta_{b}f_{0}=-2\left \langle \nabla_{b}f,\text{ }\nabla_{b}f_{0}%
\right \rangle +2V\left(  f\right)  .
\]

\end{lemma}%

%TCIMACRO{\TeXButton{Proof}{\proof} }%
%BeginExpansion
\proof
%EndExpansion
From Lemma \ref{commute}%
\[
\Delta_{b}f_{0}=\left(  \Delta_{b}f\right)  _{0}+2\underset{\alpha,\beta
=1}{\overset{n}{\sum}}\left[  \left(  A_{\alpha \beta}\varphi_{\bar{\beta}%
}\right)  _{\bar{\alpha}}+\left(  A_{\bar{\alpha}\bar{\beta}}\varphi_{\beta
}\right)  _{\alpha}\right]  .
\]
Since%
\[
\Delta_{b}f=-\left \vert \nabla_{b}f\right \vert ^{2},
\]
it follows from the commutation relation (\ref{2010a}) that
\[%
\begin{array}
[c]{lll}%
\Delta_{b}f_{0} & = & \left(  \Delta_{b}f\right)  _{0}+2\underset{\alpha
,\beta=1}{\overset{n}{\sum}}\left[  \left(  A_{\alpha \beta}f_{\bar{\beta}%
}\right)  _{\bar{\alpha}}+\left(  A_{\bar{\alpha}\bar{\beta}}f_{\beta}\right)
_{\alpha}\right] \\
& = & \left(  -\left \vert \nabla_{b}f\right \vert ^{2}\right)  _{0}%
+2\underset{\alpha,\beta=1}{\overset{n}{\sum}}\left[  \left(  A_{\alpha \beta
}f_{\bar{\beta}}\right)  _{\bar{\alpha}}+\left(  A_{\bar{\alpha}\bar{\beta}%
}f_{\beta}\right)  _{\alpha}\right] \\
& = & -2\left \langle \nabla_{b}f_{0},\text{ }\nabla_{b}f\right \rangle
+2\underset{\alpha,\beta=1}{\overset{n}{\sum}}\left[  \left(  A_{\alpha \beta
}f_{\bar{\beta}}\right)  _{\bar{\alpha}}+\left(  A_{\bar{\alpha}\bar{\beta}%
}f_{\beta}\right)  _{\alpha}+A_{\alpha \beta}f_{\bar{\alpha}}f_{\bar{\beta}%
}+A_{\bar{\alpha}\bar{\beta}}f_{\alpha}f_{\beta}\right]  .
\end{array}
\]%
%TCIMACRO{\TeXButton{End Proof}{\endproof}}%
%BeginExpansion
\endproof
%EndExpansion

\begin{lemma}
\label{V(f)} Let $(M,J,\theta)$ be a pseudohermitian $\left(  2n+1\right)
$-manifold and $u$ be a positive function with $f=\ln u.$ Suppose that%
\[
2\operatorname{Im}Qu=\left[  \Delta_{b},\text{ }T\right]  u=0.
\]
Then
\begin{equation}
V\left(  f\right)  =0. \label{116}%
\end{equation}

\end{lemma}%

%TCIMACRO{\TeXButton{Proof}{\proof} }%
%BeginExpansion
\proof
%EndExpansion
\ We compute
\begin{equation}%
\begin{array}
[c]{lll}%
V\left(  f\right)  & = & \underset{\alpha,\beta=1}{\overset{n}{\sum}}\left[
\left(  A_{\alpha \beta}f_{\bar{\beta}}\right)  _{\bar{\alpha}}+\left(
A_{\bar{\alpha}\bar{\beta}}f_{\beta}\right)  _{\alpha}+A_{\alpha \beta}%
f_{\bar{\alpha}}f_{\bar{\beta}}+A_{\bar{\alpha}\bar{\beta}}f_{\alpha}f_{\beta
}\right] \\
& = & \underset{\alpha,\beta=1}{\overset{n}{\sum}}\left[  A_{\alpha \beta
}f_{\bar{\beta}\bar{\alpha}}+A_{\alpha \beta,\bar{\alpha}}f_{\bar{\beta}%
}+A_{\bar{\alpha}\bar{\beta}}f_{\beta \alpha}+A_{\bar{\alpha}\bar{\beta}%
,\alpha}f_{\beta}+A_{\alpha \beta}f_{\bar{\alpha}}f_{\bar{\beta}}%
+A_{\bar{\alpha}\bar{\beta}}f_{\alpha}f_{\beta}\right] \\
& = & \underset{\alpha,\beta=1}{\overset{n}{\sum}}\left \{  A_{\bar{\alpha}%
\bar{\beta}}\left(  \frac{u_{\beta \alpha}}{u}-\frac{u_{\alpha}u_{\beta}}%
{u^{2}}\right)  +A_{\alpha \beta}\left(  \frac{u_{\bar{\beta}\bar{\alpha}}}%
{u}-\frac{u_{\bar{\alpha}}u_{\bar{\beta}}}{u^{2}}\right)  \right. \\
&  & \left.  +A_{\bar{\alpha}\bar{\beta},\alpha}\frac{u_{\beta}}{u}%
+A_{\alpha \beta,\bar{\alpha}}\frac{u_{\bar{\beta}}}{u}+A_{\bar{\alpha}%
\bar{\beta}}\frac{u_{\alpha}u_{\beta}}{u^{2}}+A_{\alpha \beta}\frac
{u_{\bar{\alpha}}u_{\bar{\beta}}}{u^{2}}\right \} \\
& = & \underset{\alpha,\beta=1}{\overset{n}{\sum}}\frac{1}{u}\left[  \left(
A_{\alpha \beta}u_{\bar{\beta}}\right)  _{\bar{\alpha}}+\left(  A_{\bar{\alpha
}\bar{\beta}}u_{\beta}\right)  _{\alpha}\right] \\
& = & \frac{1}{2u}\left[  \Delta_{b},\text{ }T\right]  u.
\end{array}
\label{115}%
\end{equation}

This completes the proof.
%TCIMACRO{\TeXButton{End Proof}{\endproof}}%
%BeginExpansion
\endproof
%EndExpansion

\section{CR Sub-Laplacian Comparison Theorem \ }

In this section, we give the proof of sub-Laplacian comparison theorems in a
complete noncompact pseudohermitian $(2n+1)$-manifold of vanishing
pseudohermitian torsion tensors. In order to prove Theorem \ref{t0}, we first
derive the Ricatti inequality for sub-Laplacian of Carnot-Carath\'{e}odory
distance. We refer to \cite[Corollary 3.1.]{chl} for some of computations.

\begin{lemma}
\label{la1} Let $(M,J,\theta)$ be a complete noncompact pseudohermitian
$(2n+1)$-manifold of vanishing pseudohermitian torsion with
\[
R_{\alpha \bar{\beta}}\geq k_{2}h_{\alpha \bar{\beta}}%
\]
for some constant $k_{2}$. Then, for any $x\in M$ where $r(x)$ is smooth, we have

(i) For $n=1,$
\[
\partial_{r}(\Delta_{b}r)+(\Delta_{b}r)^{2}+k_{2}\leq0.
\]

(ii) For $n\geq2,$%
\[
\partial_{r}(\Delta_{b}r)+2^{-n}(\Delta_{b}r)^{2}+k_{2}\leq0.
\]

\end{lemma}

\begin{proof}
We choose $\left \{  e_{j},e_{\widetilde{j}},T\right \}  _{j\in I_{n}}$ to be an
orthonormal frame at $q$ where $e_{\widetilde{j}}=Je_{j}$ and $e_{1}%
=\nabla_{b}r$. Since the pseudohermitian torsion is vanishing, by a result in
\cite[Corollary 2.3]{dz}, we could parallel transport such frame at $q$ to
obtain the orthonormal frame along the radial $\nabla$-geodesic $\gamma$ from
$p$ to $q$. \ Then there is an orthonormal frame $\left \{  Z_{j}%
,Z\overline{_{j}},T\right \}  _{j\in I_{n}}$ along $\gamma$ with
\[
Z_{j}=\frac{1}{\sqrt{2}}\left(  e_{j}-ie_{\widetilde{j}}\right)  .
\]
By the fact that $\gamma$ is the $\nabla$-geodesic, one can compute the
following in $Z_{1}$-direction as
\[%
\begin{array}
[c]{ccl}%
r_{11} & = & -\frac{1}{2}\left(  ie_{2}e_{1}+e_{2}e_{2}\right)  r-\left(
\nabla_{Z_{1}}Z_{1}\right)  r\\
& = & -\frac{1}{2}\left(  ie_{2}e_{1}+e_{2}e_{2}\right)  r+\frac{1}{2}\left[
i\nabla_{\left(  J\nabla_{b}r\right)  }\nabla_{b}r+J\left(  \nabla_{\left(
J\nabla_{b}r\right)  }\nabla_{b}r\right)  \right]
\end{array}
\]
and%
\[%
\begin{array}
[c]{ccl}%
r_{1\overline{1}} & = & \frac{1}{2}\left(  ie_{2}e_{1}+e_{2}e_{2}\right)
r-\left(  \nabla_{Z_{\overline{1}}}Z_{1}\right)  r\\
& = & \frac{1}{2}\left(  ie_{2}e_{1}+e_{2}e_{2}\right)  r-\frac{1}{2}\left[
i\nabla_{\left(  J\nabla_{b}r\right)  }\nabla_{b}r+J\left(  \nabla_{\left(
J\nabla_{b}r\right)  }\nabla_{b}r\right)  \right]  .
\end{array}
\]
Therefore along $\gamma$
\begin{equation}
r_{11}=-r_{1\overline{1}}. \label{302}%
\end{equation}
Moreover, by computing
\[%
\begin{array}
[c]{ccl}%
r_{1} & = & Z_{1}r\\
& = & \frac{1}{\sqrt{2}}\left(  \nabla_{b}r-iJ\nabla_{b}r\right)  r\\
& = & \frac{1}{\sqrt{2}}\left(  \left \vert \nabla_{b}r\right \vert
^{2}-i\left \langle \nabla_{b}r,J\nabla_{b}r\right \rangle \right) \\
& = & \frac{1}{\sqrt{2}}%
\end{array}
\]
and
\[%
\begin{array}
[c]{ccl}%
r_{1\overline{1}} & = & Z_{\overline{1}}Z_{1}r-\Gamma_{\overline{1}1}^{1}%
r_{1}\\
& = & Z_{\overline{1}}Z_{1}r-g_{\theta}\left(  \left[  Z_{\overline{1}}%
,Z_{1}\right]  ,Z_{1}\right)  r_{1}\\
& = & Z_{\overline{1}}Z_{1}r-\frac{1}{\sqrt{2}}g_{\theta}\left(  \left[
Z_{\overline{1}},Z_{1}\right]  ,Z_{1}\right)  ,
\end{array}
\]
we derive that $r_{1\overline{1}}$ is real by the commutation formula.
Therefore, we have
\begin{equation}
r_{0}=0 \label{304a}%
\end{equation}
along the $\nabla$-geodesic $\gamma$.

Now at the point $q$, by the facts that $r_{1}=\frac{1}{\sqrt{2}}$ and
$r_{1\overline{1}}$ is real, the equalities (\ref{302}), (\ref{304a}) and the
commutation formulas (\ref{2010a}), (\ref{2010c}), we have the following
computation as well.
\begin{equation}%
\begin{array}
[c]{ccl}%
0 & = & \frac{1}{2}\left(  \left \vert \nabla_{b}r\right \vert ^{2}\right)
_{1\overline{1}}\\
& = & \underset{\alpha}{\sum}\left(  \left \vert r_{\alpha1}\right \vert
^{2}+\left \vert r_{\alpha \overline{1}}\right \vert ^{2}+r_{\alpha1\overline{1}%
}r_{\overline{\alpha}}+r_{\overline{\alpha}1\overline{1}}r_{\alpha}\right) \\
& \geq & \left \vert r_{11}\right \vert ^{2}+\left \vert r_{1\overline{1}%
}\right \vert ^{2}+r_{11\overline{1}}r_{\overline{1}}+r_{\overline{1}%
1\overline{1}}r_{1}\\
& = & 2r_{1\overline{1}}^{2}+\left(  r_{1\overline{1}1}+ir_{10}+R_{11\overline
{1}}^{1}r_{1}\right)  r_{\overline{1}}+\left(  r_{1\overline{1}}%
-ir_{0}\right)  _{\overline{1}}r_{1}\\
& = & 2r_{1\overline{1}}^{2}+\left \langle \nabla_{b}r_{1\overline{1}}%
,\nabla_{b}r\right \rangle _{L_{\theta}}+\frac{1}{2}R_{1\overline{1}%
1\overline{1}}\\
& \geq & 2r_{1\overline{1}}^{2}+\left(  \nabla_{b}r\right)  r_{1\overline{1}%
}+\frac{1}{2}R_{1\overline{1}1\overline{1}}\\
& = & 2r_{1\overline{1}}^{2}+\left(  \nabla r\right)  r_{1\overline{1}}%
+\frac{1}{2}R_{1\overline{1}1\overline{1}}\\
& = & 2r_{1\overline{1}}^{2}+\frac{\partial r_{1\overline{1}}}{\partial
r}+\frac{1}{2}R_{1\overline{1}1\overline{1}}.
\end{array}
\label{310}%
\end{equation}

(i) For $n=1:$ Since
\[
\Delta_{b}r=r_{1\overline{1}}+r_{\overline{1}1}=2r_{1\overline{1}},
\]
it follows from (\ref{310}) that
\[
\partial_{r}(\Delta_{b}r)+(\Delta_{b}r)^{2}+R_{1\overline{1}}\leq0
\]
and then
\[
\partial_{r}(\Delta_{b}r)+(\Delta_{b}r)^{2}+k_{2}\leq0
\]
as well.

(ii) For $n\geq2:$ \ The similar computation as before, for any $j\neq1$,
\begin{equation}%
\begin{array}
[c]{ccl}%
0 & = & \frac{1}{2}\left(  \left \vert \nabla_{b}r\right \vert ^{2}\right)
_{j\overline{j}}\\
& = & \underset{\alpha}{\sum}\left(  \left \vert r_{\alpha j}\right \vert
^{2}+\left \vert r_{\alpha \overline{j}}\right \vert ^{2}+r_{\alpha j\overline
{j}}r_{\overline{\alpha}}+r_{\overline{\alpha}j\overline{j}}r_{\alpha}\right)
\\
& \geq & \left \vert r_{j\overline{j}}\right \vert ^{2}+r_{1j\overline{j}%
}r_{\overline{1}}+r_{\overline{1}j\overline{j}}r_{1}\\
& = & r_{j\overline{j}}^{2}+\left(  r_{1\overline{j}j}+ir_{10}+R_{1j\overline
{j}}^{1}r_{1}\right)  r_{\overline{1}}+r_{\overline{1}j\overline{j}}r_{1}\\
& \geq & r_{j\overline{j}}^{2}+\left \langle \nabla_{b}r_{j\overline{j}}%
,\nabla_{b}r\right \rangle _{L_{\theta}}+\frac{1}{2}R_{1\overline{1}%
j\overline{j}}\\
& = & r_{j\overline{j}}^{2}+\left(  \nabla_{b}r\right)  r_{j\overline{j}%
}+\frac{1}{2}R_{1\overline{1}j\overline{j}}\\
& = & r_{j\overline{j}}^{2}+\frac{\partial}{\partial r}r_{j\overline{j}}%
+\frac{1}{2}R_{1\overline{1}j\overline{j}}.
\end{array}
\label{315}%
\end{equation}

It follows from the inequalities (\ref{310}) and (\ref{315}) that
\begin{equation}%
\begin{array}
[c]{ccl}%
0 & \geq & \left(  2r_{1\overline{1}}^{2}+\frac{\partial r_{1\overline{1}}%
}{\partial r}+\frac{1}{2}R_{1\overline{1}1\overline{1}}\right)  +\underset
{j\neq1}{\sum}\left(  r_{j\overline{j}}^{2}+\frac{\partial}{\partial
r}r_{j\overline{j}}+\frac{1}{2}R_{1\overline{1}j\overline{j}}\right) \\
& \geq & 2^{1-n}\left(  \underset{j=1}{\overset{n}{%
%TCIMACRO{\dsum }%
%BeginExpansion
{\displaystyle \sum}
%EndExpansion
}}r_{j\overline{j}}\right)  ^{2}+\frac{\partial}{\partial r}\underset
{j=1}{\overset{n}{%
%TCIMACRO{\dsum }%
%BeginExpansion
{\displaystyle \sum}
%EndExpansion
}}r_{j\overline{j}}+\frac{1}{2}R_{1\overline{1}}.
\end{array}
\label{323}%
\end{equation}
Hence
\[
\partial_{r}(\Delta_{b}r)+2^{-n}(\Delta_{b}r)^{2}+R_{1\overline{1}}\leq0
\]
and then
\[
\partial_{r}(\Delta_{b}r)+2^{-n}(\Delta_{b}r)^{2}+k_{2}\leq0
\]
as well.
\end{proof}

Now Theorem \ref{t41}\textbf{\ }will follows from the Lemma \ref{la1} easily
(\cite{li}, \cite{w}). \ 

\begin{theorem}
\label{t41} Let $(M,J,\theta)$ be a complete pseudohermitian $(2n+1)$-manifold
of vanishing pseudohermitian torsion with
\[
R_{\alpha \bar{\beta}}\geq k_{2}h_{\alpha \bar{\beta}}%
\]
for some constant $k_{2}$. Then, for any $x\in M$ where $r(x)$ is smooth, we have

(i) $n=1$%
\begin{equation}
\Delta_{b}r\leq \left \{
\begin{array}
[c]{lc}%
\sqrt{k_{2}}\cot(\sqrt{k_{2}}r), & k_{2}>0,\\
\frac{1}{r}, & k_{2}=0,\\
\sqrt{|k_{2}|}\coth(\sqrt{|k_{2}|}r), & k_{2}<0.
\end{array}
\right.  \label{2014a}%
\end{equation}

(ii) $n\geq2:$%
\begin{equation}
\Delta_{b}r\leq \left \{
\begin{array}
[c]{lc}%
\sqrt{2^{n}k_{2}}\cot(\sqrt{2^{-n}k_{2}}), & k_{2}>0,\\
\frac{2^{n}}{r}, & k_{2}=0,\\
\sqrt{2^{n}|k_{2}|}\coth(\sqrt{2^{-n}|k_{2}|}r), & k_{2}<0.
\end{array}
\right.  \label{2014b}%
\end{equation}
Moreover, it holds on the whole manifold in the sense of distribution.
\end{theorem}

\section{CR Analogue of Yau's Gradient Estimate}

In this section, we will prove main Theorem \ref{t3} and Theorem \ref{t1}.
\ We first recall a real-valued function%
\[
F\left(  x,\text{ }t,\text{ }R,\text{ }b\right)  :M\times \left[  0,\text{
\thinspace}1\right]  \times \left(  0,\text{ }\infty \right)  \times \left(
0,\text{ }\infty \right)  \rightarrow \mathbb{R}%
\]
defined by%
\begin{equation}
F\left(  x,\text{ }t,\text{ }R,\text{ }b\right)  =t\left(  \left \vert
\nabla_{b}f\right \vert ^{2}\left(  x\right)  +bt\eta \left(  x\right)
f_{0}^{2}\left(  x\right)  \right)  , \label{2011}%
\end{equation}
where $\eta \left(  x\right)  $ $:M\rightarrow \left[  0,\text{ \thinspace
}1\right]  $ is a smooth cut-off function defined by%
\[
\eta \left(  x\right)  =\eta \left(  r\left(  x\right)  \right)  =\left \{
\begin{array}
[c]{ll}%
1, & x\in \text{ }B\left(  R\right) \\
0, & x\in M\backslash B\left(  2R\right)
\end{array}
\right.
\]
such that%
\begin{equation}
-\frac{C}{R}\eta^{\frac{1}{2}}\leq \eta^{^{\prime}}\leq0 \label{1}%
\end{equation}
and
\begin{equation}
\left \vert \eta^{^{\prime \prime}}\right \vert \leq \frac{C}{R^{2}}, \label{2}%
\end{equation}
where we denote $\frac{\partial}{\partial r}\eta$ by $\eta^{^{\prime}}$ and
$r(x)$ is the Carnot-Carath\'{e}odory distance to a fixed point $x_{0}.$

In the following calculation, the universal constant $C$ might be changed from
lines to lines.

\begin{proposition}
\label{prop1} Let $\left(  M,\text{ }J,\text{ }\theta \right)  $ be a complete
noncompact pseudohermitian $\left(  2n+1\right)  $-manifold with%
\begin{equation}
\left(  2Ric-\left(  n-2\right)  Tor\right)  \left(  Z,\text{ }Z\right)
\geq-2k\left \vert Z\right \vert ^{2} \label{113}%
\end{equation}
for all $Z\in T_{1,0}$, where $k$ is an nonnegative constant. \ Suppose that
$\left(  M,\text{ }J,\text{ }\theta \right)  $ satisfies the CR sub-Laplacian
comparison property. Then%
\[%
\begin{array}
[c]{lll}%
\Delta_{b}F & \geq & -2\left \langle \nabla_{b}f,\text{ }\nabla_{b}%
F\right \rangle +t\left[  \underset{\alpha,\beta=1}{4\overset{n}{\sum}%
}\left \vert f_{a\beta}\right \vert ^{2}+\underset{\alpha,\beta=1,\alpha
\not =\beta}{4\overset{n}{\sum}}\left \vert f_{a\bar{\beta}}\right \vert
^{2}+\frac{1}{n}\left(  \Delta_{b}f\right)  ^{2}\right. \\
&  & \left.  +\left(  n-\frac{bCt}{R}\right)  f_{0}^{2}-\left(  2k+\frac
{4}{bt\eta}\right)  \left \vert \nabla_{b}f\right \vert ^{2}-\frac{bCt}{R}%
\eta \left \vert \nabla_{b}f\right \vert ^{2}f_{0}^{2}+4bt\eta f_{0}V\left(
f\right)  \right]  .
\end{array}
\]

\end{proposition}%

%TCIMACRO{\TeXButton{Proof}{\proof} }%
%BeginExpansion
\proof
%EndExpansion
By CR sub-Laplacian comparison property,%
\[%
\begin{array}
[c]{lll}%
\Delta_{b}\eta & = & \eta^{^{\prime \prime}}+\eta^{^{\prime}}\Delta_{b}r\\
& \geq & -\frac{C}{R^{2}}-\frac{C}{R}\left(  \frac{C_{1}}{R}+C_{2}\right) \\
& \geq & -\frac{C}{R}.
\end{array}
\]

First we compute%
\begin{equation}%
\begin{array}
[c]{lll}%
\Delta_{b}\left(  bt\eta f_{0}^{2}\right)  & = & bt\left(  f_{0}^{2}\Delta
_{b}\eta+\eta \Delta_{b}f_{0}^{2}+2\left \langle \nabla_{b}\eta,\text{ }%
\nabla_{b}f_{0}^{2}\right \rangle \right) \\
& \geq & bt\left(  -\frac{C}{R}f_{0}^{2}+2\eta f_{0}\Delta_{b}f_{0}%
+2\eta \left \vert \nabla_{b}f_{0}\right \vert ^{2}+4f_{0}\left \langle \nabla
_{b}\eta,\text{ }\nabla_{b}f_{0}\right \rangle \right) \\
& \geq & bt\left(  -\frac{C}{R}f_{0}^{2}+2\eta f_{0}\Delta_{b}f_{0}%
+2\eta \left \vert \nabla_{b}f_{0}\right \vert ^{2}-4\left \vert f_{0}\right \vert
\left \vert \nabla_{b}\eta \right \vert \left \vert \nabla_{b}f_{0}\right \vert
\right) \\
& \geq & bt[-\frac{C}{R}f_{0}^{2}+2\eta f_{0}\Delta_{b}f_{0}+\left(
2\eta-2\cdot \frac{1}{2}\eta \right)  \left \vert \nabla_{b}f_{0}\right \vert
^{2}-2\cdot2\eta^{-1}\left \vert \nabla_{b}\eta \right \vert ^{2}f_{0}^{2}]\\
& \geq & bt[-\frac{C}{R}f_{0}^{2}+2\eta f_{0}\Delta_{b}f_{0}+\left(
2\eta-2\cdot \frac{1}{2}\eta \right)  \left \vert \nabla_{b}f_{0}\right \vert
^{2}],
\end{array}
\label{3}%
\end{equation}
where we use the Young's inequality and the inequality (\ref{1}) which implies
that%
\[
\eta^{-1}\left \vert \nabla_{b}\eta \right \vert ^{2}\leq \frac{C}{R^{2}}.
\]

Second, it follows from assumption (\ref{113}), Lemma \ref{Bochner inequality}
and (\ref{3}) that
\[%
\begin{array}
[c]{lll}%
\Delta_{b}F & = & t\left(  \Delta_{b}\left \vert \nabla_{b}f\right \vert
^{2}+\Delta_{b}\left(  bt\eta f_{0}^{2}\right)  \right) \\
& \geq & t\left(  \underset{\alpha,\beta=1}{4\overset{n}{\sum}}\left \vert
f_{a\beta}\right \vert ^{2}+\underset{\alpha,\beta=1,\alpha \not =\beta
}{4\overset{n}{\sum}}\left \vert f_{a\bar{\beta}}\right \vert ^{2}+\frac{1}%
{n}\left(  \Delta_{b}f\right)  ^{2}+nf_{0}^{2}+2\left \langle \nabla
_{b}f,\text{ }\nabla_{b}\Delta_{b}f\right \rangle \right. \\
&  & \left.  -2\left(  k+\frac{1}{\nu}\right)  \left \vert \nabla
_{b}f\right \vert ^{2}-2\nu \left \vert \nabla_{b}f_{0}\right \vert ^{2}%
-\frac{bCt}{R}f_{0}^{2}+2bt\eta f_{0}\Delta_{b}f_{0}+2\cdot \frac{bt}{2}%
\eta \left \vert \nabla_{b}f_{0}\right \vert ^{2}\right) \\
& \geq & t\left[  \underset{\alpha,\beta=1}{4\overset{n}{\sum}}\left \vert
f_{a\beta}\right \vert ^{2}+\underset{\alpha,\beta=1,\alpha \not =\beta
}{4\overset{n}{\sum}}\left \vert f_{a\bar{\beta}}\right \vert ^{2}+\frac{1}%
{n}\left(  \Delta_{b}f\right)  ^{2}+\left(  n-\frac{bCt}{R}\right)  f_{0}%
^{2}+2\left \langle \nabla_{b}f,\text{ }\nabla_{b}\Delta_{b}f\right \rangle
\right. \\
&  & \left.  -2\left(  k+\frac{1}{\nu}\right)  \left \vert \nabla
_{b}f\right \vert ^{2}+2\left(  \frac{bt}{2}\eta-\nu \right)  \left \vert
\nabla_{b}f_{0}\right \vert ^{2}+2bt\eta f_{0}\Delta_{b}f_{0}\right]  .
\end{array}
\]

Then taking $\nu=\frac{bt\eta}{2}$,%
\begin{equation}%
\begin{array}
[c]{lll}%
\Delta_{b}F & \geq & t\left[  \underset{\alpha,\beta=1}{4\overset{n}{\sum}%
}\left \vert f_{a\beta}\right \vert ^{2}+\underset{\alpha,\beta=1,\alpha
\not =\beta}{4\overset{n}{\sum}}\left \vert f_{a\bar{\beta}}\right \vert
^{2}+\frac{1}{n}\left(  \Delta_{b}f\right)  ^{2}+\left(  n-\frac{bCt}%
{R}\right)  f_{0}^{2}\right. \\
&  & \left.  -2\left(  k+\frac{2}{bt\eta}\right)  \left \vert \nabla
_{b}f\right \vert ^{2}+2\left \langle \nabla_{b}f,\text{ }\nabla_{b}\Delta
_{b}f\right \rangle +2bt\eta f_{0}\Delta_{b}f_{0}\right]  .
\end{array}
\label{4}%
\end{equation}

Finally, by Lemma \ref{commute f}%
\[%
\begin{array}
[c]{ll}
& 2\left \langle \nabla_{b}f,\text{ }\nabla_{b}\Delta_{b}f\right \rangle
+2bt\eta f_{0}\Delta_{b}f_{0}\\
= & 2\left \langle \nabla_{b}f,\text{ }\nabla_{b}\left(  -\left \vert \nabla
_{b}f\right \vert ^{2}\right)  \right \rangle +2bt\eta f_{0}\left[
-2\left \langle \nabla_{b}f,\text{ }\nabla_{b}f_{0}\right \rangle +2V\left(
f\right)  \right] \\
= & -2\left \langle \nabla_{b}f,\text{ }\nabla_{b}\left(  \frac{F}{t}-bt\eta
f_{0}^{2}\right)  \right \rangle -4bt\eta f_{0}\left \langle \nabla_{b}f,\text{
}\nabla_{b}f_{0}\right \rangle +4bt\eta f_{0}V\left(  f\right) \\
= & -\frac{2}{t}\left \langle \nabla_{b}f,\text{ }\nabla_{b}F\right \rangle
+2bt\left \langle \nabla_{b}f,\text{ }\nabla_{b}\left(  \eta f_{0}^{2}\right)
\right \rangle -4bt\eta f_{0}\left \langle \nabla_{b}f,\text{ }\nabla_{b}%
f_{0}\right \rangle +4bt\eta f_{0}V\left(  f\right) \\
= & -\frac{2}{t}\left \langle \nabla_{b}f,\text{ }\nabla_{b}F\right \rangle
+2btf_{0}^{2}\left \langle \nabla_{b}f,\text{ }\nabla_{b}\eta \right \rangle
+4bt\eta f_{0}V\left(  f\right)
\end{array}
\]

Now by Young's inequality, we have%
\begin{equation}%
\begin{array}
[c]{ll}
& 2\left \langle \nabla_{b}f,\text{ }\nabla_{b}\Delta_{b}f\right \rangle
+2bt\eta f_{0}\Delta_{b}f_{0}\\
= & -\frac{2}{t}\left \langle \nabla_{b}f,\text{ }\nabla_{b}F\right \rangle
+2btf_{0}^{2}\left \langle \nabla_{b}f,\text{ }\nabla_{b}\eta \right \rangle
+4bt\eta f_{0}V\left(  f\right) \\
\geq & -\frac{2}{t}\left \langle \nabla_{b}f,\text{ }\nabla F\right \rangle
-2btf_{0}^{2}\left \vert \nabla_{b}f\right \vert \left \vert \nabla_{b}%
\eta \right \vert +4bt\eta f_{0}V\left(  f\right) \\
\geq & -\frac{2}{t}\left \langle \nabla_{b}f,\text{ }\nabla F\right \rangle
-\frac{2Cbt}{R}f_{0}^{2}\left \vert \nabla_{b}f\right \vert \eta^{\frac{1}{2}%
}+4bt\eta f_{0}V\left(  f\right) \\
\geq & -\frac{2}{t}\left \langle \nabla_{b}f,\text{ }\nabla F\right \rangle
-\frac{Cbt}{R}f_{0}^{2}-\frac{Cbt}{R}\eta f_{0}^{2}\left \vert \nabla
_{b}f\right \vert ^{2}+4bt\eta f_{0}V\left(  f\right)  .
\end{array}
\label{5}%
\end{equation}

Substituting (\ref{5}) into (\ref{4}),%
\[%
\begin{array}
[c]{lll}%
\Delta_{b}F & \geq & -2\left \langle \nabla_{b}f,\text{ }\nabla F\right \rangle
+t\left[  \underset{\alpha,\beta=1}{4\overset{n}{\sum}}\left \vert f_{a\beta
}\right \vert ^{2}+\underset{\alpha,\beta=1,\alpha \not =\beta}{4\overset
{n}{\sum}}\left \vert f_{a\bar{\beta}}\right \vert ^{2}+\frac{1}{n}\left(
\Delta_{b}f\right)  ^{2}\right. \\
&  & \left.  +\left(  n-\frac{bCt}{R}\right)  f_{0}^{2}-2\left(  k+\frac
{2}{bt\eta}\right)  \left \vert \nabla_{b}f\right \vert ^{2}-\frac{Cbt}{R}\eta
f_{0}^{2}\left \vert \nabla_{b}f\right \vert ^{2}+4bt\eta f_{0}V\left(
f\right)  \right]  .
\end{array}
\]%
%TCIMACRO{\TeXButton{End Proof}{\endproof}}%
%BeginExpansion
\endproof
%EndExpansion

\begin{proposition}
\label{prop2} \ Let $\left(  M,\text{ }J,\text{ }\theta \right)  $ be a
complete noncompact pseudohermitian $\left(  2n+1\right)  $-manifold with%
\[
\left(  2Ric-\left(  n-2\right)  Tor\right)  \left(  Z,\text{ }Z\right)
\geq-2k\left \vert Z\right \vert ^{2}%
\]
for all $Z\in T_{1,0}$, where $k$ is an nonnegative constant. \ Suppose that
$\left(  M,\text{ }J,\text{ }\theta \right)  $ satisfies the CR sub-Laplacian
comparison property. Then for all $a\neq0$%
\begin{equation}%
\begin{array}
[c]{lll}%
t\eta \Delta_{b}\left(  \eta F\right)  & \geq & \frac{1}{na^{2}}\left(  \eta
F\right)  ^{2}-\frac{C}{R}\left(  \eta F\right)  +2t\eta \left \langle
\nabla_{b}\eta,\text{ }\nabla_{b}F\right \rangle -2t\eta^{2}\left \langle
\nabla_{b}f,\text{ }\nabla_{b}F\right \rangle \\
&  & +4t^{2}\eta^{2}\left(  \underset{\alpha,\beta=1}{\overset{n}{\sum}%
}\left \vert f_{a\beta}\right \vert ^{2}+\underset{\alpha,\beta=1,\alpha
\not =\beta}{\overset{n}{\sum}}\left \vert f_{a\bar{\beta}}\right \vert
^{2}\right) \\
&  & +\left[  n-\frac{bC}{R}-\left(  \frac{2b}{na^{2}}+\frac{bC}{R}\right)
\left(  \eta F\right)  \right]  t^{2}\eta^{2}f_{0}^{2}\\
&  & +\left[  -\frac{2\left(  1+a\right)  }{na^{2}}\left(  \eta F\right)
-2k-\frac{4}{b}\right]  t\eta \left \vert \nabla_{b}f\right \vert ^{2}%
+4bt^{3}\eta^{3}f_{0}V\left(  f\right)  .
\end{array}
\label{7}%
\end{equation}

\end{proposition}%

%TCIMACRO{\TeXButton{Proof}{\proof} }%
%BeginExpansion
\proof
%EndExpansion
\  \ By using Proposition \ref{prop1}, we first compute
\[%
\begin{array}
[c]{lll}%
\Delta_{b}\left(  \eta F\right)  & = & \left(  \Delta_{b}\eta \right)
F+2\left \langle \nabla_{b}\eta,\text{ }\nabla_{b}F\right \rangle +\eta
\Delta_{b}F\\
& \geq & -\frac{C}{R}F+2\left \langle \nabla_{b}\eta,\text{ }\nabla
_{b}F\right \rangle -2\eta \left \langle \nabla_{b}f,\text{ }\nabla
F\right \rangle \\
&  & +t\eta \left[  4\left(  \underset{\alpha,\beta=1}{\overset{n}{\sum}%
}\left \vert f_{a\beta}\right \vert ^{2}+\underset{\alpha,\beta=1,\alpha
\not =\beta}{\overset{n}{\sum}}\left \vert f_{a\bar{\beta}}\right \vert
^{2}\right)  +\frac{1}{n}\left(  \Delta_{b}f\right)  ^{2}+\left(  n-\frac
{bCt}{R}\right)  f_{0}^{2}\right. \\
&  & \left.  -2\left(  k+\frac{2}{bt\eta}\right)  \left \vert \nabla
_{b}f\right \vert ^{2}-\frac{Cbt}{R}\eta f_{0}^{2}\left \vert \nabla
_{b}f\right \vert ^{2}+4bt\eta f_{0}V\left(  f\right)  \right]
\end{array}
\]
and for each $a\neq0$%
\[%
\begin{array}
[c]{lll}%
\left(  \Delta_{b}f\right)  ^{2} & = & \left(  -\left \vert \nabla
_{b}f\right \vert ^{2}\right)  ^{2}\\
& = & \left(  \frac{1}{at}F-\frac{1}{a}\left \vert \nabla_{b}f\right \vert
^{2}-\frac{1}{a}bt\eta f_{0}^{2}-\left \vert \nabla_{b}f\right \vert
^{2}\right)  ^{2}\\
& = & \left(  \frac{1}{at}F-\frac{a+1}{a}\left \vert \nabla_{b}f\right \vert
^{2}-\frac{1}{a}bt\eta f_{0}^{2}\right)  ^{2}\\
& = & \frac{1}{a^{2}t^{2}}F^{2}+\left(  \frac{a+1}{a}\right)  ^{2}\left \vert
\nabla_{b}f\right \vert ^{4}+\frac{1}{a^{2}}b^{2}t^{2}\eta^{2}f_{0}^{4}\\
&  & -\frac{2\left(  a+1\right)  }{a^{2}t}F\left \vert \nabla_{b}f\right \vert
^{2}-\frac{2b}{a^{2}}\eta Ff_{0}^{2}+\frac{2\left(  a+1\right)  bt}{a^{2}}%
\eta \left \vert \nabla_{b}f\right \vert ^{2}f_{0}^{2}\\
& \geq & \frac{1}{a^{2}t^{2}}F^{2}-\frac{2\left(  a+1\right)  }{a^{2}%
t}F\left \vert \nabla_{b}f\right \vert ^{2}-\frac{2b}{a^{2}}\eta Ff_{0}^{2}.
\end{array}
\]

Then%
\[%
\begin{array}
[c]{lll}%
\Delta_{b}\left(  \eta F\right)  & \geq & \frac{1}{na^{2}t}\eta F^{2}-\frac
{C}{R}F+2\left \langle \nabla_{b}\eta,\text{ }\nabla_{b}F\right \rangle
-2\eta \left \langle \nabla_{b}f,\text{ }\nabla F\right \rangle \\
&  & +4t\eta \left(  \underset{\alpha,\beta=1}{\overset{n}{\sum}}\left \vert
f_{a\beta}\right \vert ^{2}+\underset{\alpha,\beta=1,\alpha \not =\beta
}{\overset{n}{\sum}}\left \vert f_{a\bar{\beta}}\right \vert ^{2}\right) \\
&  & +\left(  n-\frac{bCt}{R}-\frac{2b}{na^{2}}\eta F\right)  t\eta f_{0}%
^{2}+\left(  -\frac{2\left(  1+a\right)  }{na^{2}}\eta F-2kt\eta-\frac{4}%
{b}\right)  \left \vert \nabla_{b}f\right \vert ^{2}\\
&  & -\frac{Cb}{R}\left(  t\eta \left \vert \nabla_{b}f\right \vert ^{2}\right)
\left(  t\eta f_{0}^{2}\right)  +4bt^{2}\eta^{2}f_{0}V\left(  f\right)  .
\end{array}
\]
Hence%
\begin{equation}%
\begin{array}
[c]{lll}%
\Delta_{b}\left(  \eta F\right)  & \geq & \frac{1}{na^{2}t}\eta F^{2}-\frac
{C}{R}F+2\left \langle \nabla_{b}\eta,\text{ }\nabla_{b}F\right \rangle
-2\eta \left \langle \nabla_{b}f,\text{ }\nabla F\right \rangle \\
&  & +4t\eta \left(  \underset{\alpha,\beta=1}{\overset{n}{\sum}}\left \vert
f_{a\beta}\right \vert ^{2}+\underset{\alpha,\beta=1,\alpha \not =\beta
}{\overset{n}{\sum}}\left \vert f_{a\bar{\beta}}\right \vert ^{2}\right) \\
&  & +\left[  n-\frac{bCt}{R}-\left(  \frac{2b}{na^{2}}+\frac{bC}{R}\right)
\eta F\right]  t\eta f_{0}^{2}\\
&  & +\left(  -\frac{2\left(  1+a\right)  }{na^{2}}\eta F-2kt\eta-\frac{4}%
{b}\right)  \left \vert \nabla_{b}f\right \vert ^{2}+4bt^{2}\eta^{2}%
f_{0}V\left(  f\right)  .
\end{array}
\label{6}%
\end{equation}

Finally, multiply $t\eta \,$\ on the both sides of (\ref{6}) and note that
$t\leq1,$ $\eta \leq1$%
\[%
\begin{array}
[c]{lll}%
t\eta \Delta_{b}\left(  \eta F\right)  & \geq & \frac{1}{na^{2}}\left(  \eta
F\right)  ^{2}-\frac{C}{R}\eta F+2t\eta \left \langle \nabla_{b}\eta,\text{
}\nabla_{b}F\right \rangle -2t\eta^{2}\left \langle \nabla_{b}f,\text{ }%
\nabla_{b}F\right \rangle \\
&  & +4t^{2}\eta^{2}\left(  \underset{\alpha,\beta=1}{\overset{n}{\sum}%
}\left \vert f_{a\beta}\right \vert ^{2}+\underset{\alpha,\beta=1,\alpha
\not =\beta}{\overset{n}{\sum}}\left \vert f_{a\bar{\beta}}\right \vert
^{2}\right) \\
&  & +\left[  n-\frac{bC}{R}-\left(  \frac{2b}{na^{2}}+\frac{bC}{R}\right)
\left(  \eta F\right)  \right]  t^{2}\eta^{2}f_{0}^{2}\\
&  & +\left[  -\frac{2\left(  1+a\right)  }{na^{2}}\left(  \eta F\right)
-2k-\frac{4}{b}\right]  t\eta \left \vert \nabla_{b}f\right \vert ^{2}%
+4bt^{3}\eta^{3}f_{0}V\left(  f\right)  .
\end{array}
\]%
%TCIMACRO{\TeXButton{End Proof}{\endproof}}%
%BeginExpansion
\endproof
%EndExpansion

\begin{proposition}
\label{prop3} Let $\left(  M,\text{ }J,\text{ }\theta \right)  $ be a complete
noncompact pseudohermitian $\left(  2n+1\right)  $-manifold with%
\[
\left(  2Ric-\left(  n-2\right)  Tor\right)  \left(  Z,\text{ }Z\right)
\geq-2k\left \vert Z\right \vert ^{2}%
\]
for all $Z\in T_{1,0}$, where $k$ is an nonnegative constant. \ Suppose that
$\left(  M,\text{ }J,\text{ }\theta \right)  $ satisfies the CR sub-Laplacian
comparison property. Let $b$, $R$ be fixed, and $p\left(  t\right)  \in
B\left(  2R\right)  $ be the maximal point of $\eta F$ for each $t\in \left(
0,1\right]  $. Then at $\left(  \text{ }p\left(  t\right)  ,t\right)  $ we
have%
\begin{equation}%
\begin{array}
[c]{lll}%
0 & \geq & \left(  \frac{1}{na^{2}}-\frac{C}{R}\right)  \left(  \eta F\right)
^{2}-\frac{3C}{R}\left(  \eta F\right)  +4t^{2}\eta^{2}\left(  \underset
{\alpha,\beta=1}{\overset{n}{\sum}}\left \vert f_{a\beta}\right \vert
^{2}+\underset{\alpha,\beta=1,\alpha \not =\beta}{\overset{n}{\sum}}\left \vert
f_{a\bar{\beta}}\right \vert ^{2}\right) \\
&  & +\left[  n-\frac{bC}{R}-\left(  \frac{2b}{na^{2}}+\frac{bC}{R}\right)
\left(  \eta F\right)  \right]  t^{2}\eta^{2}f_{0}^{2}\\
&  & +\left[  -\frac{2\left(  1+a\right)  }{na^{2}}\left(  \eta F\right)
-2k-\frac{4}{b}-\frac{C}{R}\right]  t\eta \left \vert \nabla_{b}f\right \vert
^{2}+4bt^{3}\eta^{3}f_{0}V\left(  f\right)  .
\end{array}
\label{10}%
\end{equation}

\end{proposition}%

%TCIMACRO{\TeXButton{Proof}{\proof} }%
%BeginExpansion
\proof
%EndExpansion
\ Since $\left(  \eta F\right)  \left(  p\left(  t\right)  ,\text{ }t,\text{
}R,\text{ }b\right)  =\underset{x\in B\left(  2R\right)  }{\max}\left(  \eta
F\right)  \left(  x,\text{ }t,\text{ }R,\text{ }b\right)  $, at a critical
point $\left(  p\left(  t\right)  ,\text{ }t\right)  $ of $\left(  \eta
F\right)  \left(  x,\text{ }t,\text{ }R,\text{ }b\right)  $, we have%
\[
\nabla_{b}\left(  \eta F\right)  \left(  p\left(  t\right)  ,\text{ }t,\text{
}R,\text{ }b\right)  =0.
\]
This implies that%
\begin{equation}
F\nabla_{b}\eta+\eta \nabla_{b}F=0 \label{critical point}%
\end{equation}
at $\left(  p\left(  t\right)  ,\text{ }t\right)  .$ On the other hand,%
\begin{equation}
\Delta_{b}\left(  \eta F\right)  \left(  p\left(  t\right)  ,\text{ }t,\text{
}R,\text{ }b\right)  \leq0 \label{max-principle}%
\end{equation}
at $\left(  p\left(  t\right)  ,\text{ }t\right)  .$

Now we apply (\ref{critical point}) to $2t\eta \left \langle \nabla_{b}%
\eta,\text{ }\nabla_{b}F\right \rangle $ and $-2t\eta^{2}\left \langle
\nabla_{b}f,\text{ }\nabla_{b}F\right \rangle $ in (\ref{7}), we can derive the
following estimates.
\begin{equation}%
\begin{array}
[c]{lll}%
2t\eta \left \langle \nabla_{b}\eta,\text{ }\nabla_{b}F\right \rangle  & = &
-2tF\left \vert \nabla_{b}\eta \right \vert ^{2}\\
& \geq & -\frac{2tC}{R^{2}}\eta F\\
& \geq & -\frac{2C}{R}\eta F
\end{array}
\label{8}%
\end{equation}
and
\begin{equation}%
\begin{array}
[c]{lll}%
-2t\eta^{2}\left \langle \nabla_{b}f,\text{ }\nabla_{b}F\right \rangle  & = &
2t\eta F\left \langle \nabla_{b}f,\text{ }\nabla_{b}\eta \right \rangle \\
& \geq & -2t\left(  \eta F\right)  \left \vert \nabla_{b}f\right \vert
\left \vert \text{ }\nabla_{b}\eta \right \vert \\
& \geq & -\frac{2tC}{R}\left(  \eta F\right)  \eta^{\frac{1}{2}}\left \vert
\nabla_{b}f\right \vert \\
& \geq & -\frac{Ct}{R}\left(  \eta F\right)  ^{2}-\frac{C}{R}t\eta \left \vert
\nabla_{b}f\right \vert ^{2}.
\end{array}
\label{9}%
\end{equation}
Here we have applied the Young's inequality for (\ref{9}).

Finally, substituting (\ref{max-principle}), (\ref{8}) and (\ref{9}) into
(\ref{7}) in Proposition \ref{prop2}, and noting that $t\leq1$,%
\[%
\begin{array}
[c]{lll}%
0 & \geq & \left(  \frac{1}{na^{2}}-\frac{C}{R}\right)  \left(  \eta F\right)
^{2}-\frac{3C}{R}\left(  \eta F\right)  +4t^{2}\eta^{2}\left(  \underset
{\alpha,\beta=1}{\overset{n}{\sum}}\left \vert f_{a\beta}\right \vert
^{2}+\underset{\alpha,\beta=1,\alpha \not =\beta}{\overset{n}{\sum}}\left \vert
f_{a\bar{\beta}}\right \vert ^{2}\right) \\
&  & +\left[  n-\frac{bC}{R}-\left(  \frac{2b}{na^{2}}+\frac{bC}{R}\right)
\left(  \eta F\right)  \right]  t^{2}\eta^{2}f_{0}^{2}\\
&  & +\left[  -\frac{2\left(  1+a\right)  }{na^{2}}\left(  \eta F\right)
-2k-\frac{4}{b}-\frac{C}{R}\right]  t\eta \left \vert \nabla_{b}f\right \vert
^{2}+4bt^{3}\eta^{3}f_{0}V\left(  f\right)  .
\end{array}
\]
This completes the proof.
%TCIMACRO{\TeXButton{End Proof}{\endproof}}%
%BeginExpansion
\endproof
%EndExpansion

Now, we are ready to prove our main theorems.

\textbf{Proof of Theorem\  \ref{t1} }:%

%TCIMACRO{\TeXButton{Proof}{\proof} }%
%BeginExpansion
\proof
%EndExpansion
We observe that
\begin{equation}
V\left(  f\right)  =0 \label{2011b}%
\end{equation}
by assumption (\ref{111}) and Lemma \ref{V(f)}.

We begin by substituting (\ref{2011b}) into (\ref{10}) in Proposition
\ref{prop3} at the maximum point $\left(  p(t),t\right)  $. \ Hence
\begin{equation}%
\begin{array}
[c]{lll}%
0 & \geq & \left(  \frac{1}{na^{2}}-\frac{C}{R}\right)  \left[  \left(  \eta
F\right)  \right]  ^{2}-\frac{3C}{R}\left[  \left(  \eta F\right)  \right] \\
&  & +\left[  n-\frac{bC}{R}-\left(  \frac{2b}{na^{2}}+\frac{bC}{R}\right)
\left(  \eta F\right)  \right]  t^{2}\eta^{2}f_{0}^{2}\\
&  & +\left[  -\frac{2\left(  1+a\right)  }{na^{2}}\left(  \eta F\right)
-2k-\frac{4}{b}-\frac{C}{R}\right]  t\eta \left \vert \nabla_{b}f\right \vert
^{2}\\
&  & +4t_{0}^{2}\eta^{2}\left(  \underset{\alpha,\beta=1}{\overset{n}{\sum}%
}\left \vert f_{a\beta}\right \vert ^{2}+\underset{\alpha,\beta=1,\alpha
\not =\beta}{\overset{n}{\sum}}\left \vert f_{a\bar{\beta}}\right \vert
^{2}\right)  .
\end{array}
\label{contradiction}%
\end{equation}

We claim at $t=1$%
\begin{equation}
\left(  \eta F\right)  \left(  p\left(  1\right)  ,\text{ }1,\text{ }R,\text{
}b\right)  <\frac{na^{2}}{-2\left(  1+a\right)  }\left(  2k+\frac{4}{b}%
+\frac{C}{R}\right)  \label{a}%
\end{equation}
for a large enough $R$ which to be determined later. Here $\left(  1+a\right)
<0$ for some $a$ to be chosen later (say $1+a=-\frac{5+2bk}{n})$.

We prove it by contradiction. Suppose not, that is
\[
\left(  \eta F\right)  \left(  p\left(  1\right)  ,\text{ }1,\text{ }R,\text{
}b\right)  \geq \frac{na^{2}}{-2\left(  1+a\right)  }\left(  2k+\frac{4}%
{b}+\frac{C}{R}\right)  .
\]
Since $\left(  \eta F\right)  \left(  p\left(  t\right)  ,\text{ }t,\text{
}R,\text{ }b\right)  $ is continuous in the variable $t$ and $\left(  \eta
F\right)  \left(  p\left(  0\right)  ,\text{ }0,\text{ }R,\text{ }b\right)
=0$, by Intermediate-value theorem there exists a $t_{0}\in \left(  0,\text{
}1\right]  $ such that%
\begin{equation}
\left(  \eta F\right)  \left(  p\left(  t_{0}\right)  ,\text{ }t_{0},\text{
}R,\text{ }b\right)  =\frac{na^{2}}{-2\left(  1+a\right)  }\left(  2k+\frac
{4}{b}+\frac{C}{R}\right)  . \label{1122}%
\end{equation}

Now we apply (\ref{contradiction}) at the point $\left(  p\left(
t_{0}\right)  ,t_{0}\right)  $, denoted by $\left(  p_{0},t_{0}\right)  $.
\ We have by using (\ref{1122})
\begin{equation}%
\begin{array}
[c]{lll}%
0 & \geq & \left(  \frac{1}{na^{2}}-\frac{C}{R}\right)  \left[  \left(  \eta
F\right)  \left(  p_{0},t_{0}\right)  \right]  ^{2}-\frac{3C}{R}\left[
\left(  \eta F\right)  \left(  p_{0},t_{0}\right)  \right] \\
&  & +\left[  n-\frac{bC}{R}-\left(  \frac{2b}{na^{2}}+\frac{bC}{R}\right)
\left(  \eta F\right)  \left(  p_{0},t_{0}\right)  \right]  t^{2}\eta^{2}%
f_{0}^{2}\\
&  & +4t_{0}^{2}\eta^{2}\left(  \underset{\alpha,\beta=1}{\overset{n}{\sum}%
}\left \vert f_{a\beta}\right \vert ^{2}+\underset{\alpha,\beta=1,\alpha
\not =\beta}{\overset{n}{\sum}}\left \vert f_{a\bar{\beta}}\right \vert
^{2}\right)  .
\end{array}
\label{21}%
\end{equation}
\ Moreover, \ we compute%
\begin{equation}%
\begin{array}
[c]{l}%
\left[  \left(  \frac{1}{na^{2}}-\frac{C}{R}\right)  \left(  \eta F\right)
\left(  p_{0},t_{0}\right)  -\frac{3C}{R}\right] \\
=\left[  \left(  \frac{1}{na^{2}}-\frac{C}{R}\right)  \left(  \frac{na^{2}%
}{-2\left(  1+a\right)  }\right)  \left(  2k+\frac{4}{b}+\frac{C}{R}\right)
-\frac{3C}{R}\right] \\
=\left \{  \frac{-1}{2(1+a)}\left(  2k+\frac{4}{b}\right)  -\frac{C}{R}\left[
\frac{na^{2}}{-2\left(  1+a\right)  }\left(  2k+\frac{4}{b}+\frac{C}%
{R}\right)  +\frac{1}{2\left(  1+a\right)  }+3\right]  \right \}
\end{array}
\label{12}%
\end{equation}
and
\begin{equation}%
\begin{array}
[c]{l}%
\left[  n-\frac{bC}{R}-\left(  \frac{2b}{na^{2}}+\frac{bC}{R}\right)  \left(
\eta F\right)  \left(  p_{0},t_{0}\right)  \right] \\
=n-\frac{bC}{R}-\left(  \frac{2b}{na^{2}}+\frac{bC}{R}\right)  \left(
\frac{na^{2}}{-2\left(  1+a\right)  }\right)  \left(  2k+\frac{4}{b}+\frac
{C}{R}\right) \\
=n-\frac{bC}{R}+\frac{b}{\left(  1+a\right)  }\left(  2k+\frac{4}{b}+\frac
{C}{R}\right)  +\frac{bC}{R}\left(  \frac{na^{2}}{2\left(  1+a\right)
}\right)  \left(  2k+\frac{4}{b}+\frac{C}{R}\right) \\
=\left(  n+\frac{4}{1+a}+\frac{2bk}{1+a}\right)  +\frac{C}{R}\left[
-\frac{ab}{1+a}+\frac{na^{2}b}{2(1+a)}\left(  2k+\frac{4}{b}+\frac{C}%
{R}\right)  \right]  .
\end{array}
\label{13}%
\end{equation}

Now we choose $a$ such that
\[
(1+a)<-\frac{4+2bk}{n}%
\]
and then
\[
\left(  n+\frac{4}{1+a}+\frac{2bk}{1+a}\right)  >0.
\]
In particular, we let
\begin{equation}
1+a=-\frac{5+2bk}{n}. \label{b}%
\end{equation}
Then for $R=R(b,k)$ large enough, one obtains
\[
\left[  \left(  \frac{1}{na^{2}}-\frac{C}{R}\right)  \left(  \eta F\right)
\left(  p_{0},t_{0}\right)  -\frac{3C}{R}\right]  >0
\]
and
\[
\left[  n-\frac{bC}{R}-\left(  \frac{2b}{na^{2}}+\frac{bC}{R}\right)  \left(
\eta F\right)  \left(  p_{0},t_{0}\right)  \right]  >0.
\]
This leads to a contradiction with (\ref{21}). Hence from (\ref{a}) and
(\ref{b})%
\[
\left(  \eta F\right)  \left(  1,\text{ }p\left(  1\right)  ,\text{ }R,\text{
}b\right)  <\frac{\left(  n+5+2bk\right)  ^{2}}{2\left(  5+2bk\right)
}\left(  2k+\frac{4}{b}+\frac{C}{R}\right)  .
\]
This implies
\[
\underset{x\in B\left(  2R\right)  }{\max}\left(  \left \vert \nabla
_{b}f\right \vert ^{2}+b\eta f_{0}^{2}\right)  \left(  x\right)  <\frac{\left(
n+5+2bk\right)  ^{2}}{2\left(  5+2bk\right)  }\left(  2k+\frac{4}{b}+\frac
{C}{R}\right)  .
\]
When we fix on the set $x\in B\left(  R\right)  $, we obtain%
\[
\left \vert \nabla_{b}f\right \vert ^{2}+bf_{0}^{2}<\frac{\left(
n+5+2bk\right)  ^{2}}{2\left(  5+2bk\right)  }\left(  2k+\frac{4}{b}+\frac
{C}{R}\right)
\]
on $B\left(  R\right)  $.

This completes the proof.\
%TCIMACRO{\TeXButton{endproof}{\hfill}}%
%BeginExpansion
\hfill
%EndExpansion

\bigskip

Next we prove Theorem \ref{t3}. The proof is similar to Theorem \ref{t1}. \ 

\textbf{Proof of Theorem \ref{t3} }:%

%TCIMACRO{\TeXButton{Proof}{\proof} }%
%BeginExpansion
\proof
%EndExpansion
Firstly, we recall (Proposition \ref{prop2}) that%
\begin{equation}%
\begin{array}
[c]{ll}
& t\eta \Delta_{b}\left(  \eta F\right) \\
\geq & \frac{1}{na^{2}}\left(  \eta F\right)  ^{2}-\frac{C}{R}\left(  \eta
F\right)  +2t\eta \left \langle \nabla_{b}\eta,\text{ }\nabla_{b}F\right \rangle
-2t\eta^{2}\left \langle \nabla_{b}f,\text{ }\nabla_{b}F\right \rangle \\
& +4t^{2}\eta^{2}\left(  \underset{\alpha,\beta=1}{\overset{n}{\sum}%
}\left \vert f_{a\beta}\right \vert ^{2}+\underset{\alpha,\beta=1,\alpha
\not =\beta}{\overset{n}{\sum}}\left \vert f_{a\bar{\beta}}\right \vert
^{2}\right) \\
& +\left[  n-\frac{bC}{R}-\left(  \frac{2b}{na^{2}}+\frac{bC}{R}\right)
\left(  \eta F\right)  \right]  t^{2}\eta^{2}f_{0}^{2}\\
& +\left[  -\frac{2\left(  1+a\right)  }{na^{2}}\left(  \eta F\right)
-2k-\frac{4}{b}\right]  t\eta \left \vert \nabla_{b}f\right \vert ^{2}%
+4bt^{3}\eta^{3}f_{0}V\left(  f\right)  .
\end{array}
\label{15}%
\end{equation}

Now we need to deal with the term $4bt^{3}\eta^{3}f_{0}V\left(  f\right)  $ in
(\ref{15}).%
\begin{equation}%
\begin{array}
[c]{l}%
4bt^{3}\eta^{3}f_{0}V\left(  f\right) \\
=4bt^{3}\eta^{3}f_{0}\underset{\alpha,\beta=1}{\overset{n}{\sum}}\left[
\left(  A_{\alpha \beta}f_{\bar{\beta}}\right)  _{\bar{\alpha}}+\left(
A_{\bar{\alpha}\bar{\beta}}f_{\beta}\right)  _{\alpha}+A_{\alpha \beta}%
f_{\bar{\alpha}}f_{\bar{\beta}}+A_{\bar{\alpha}\bar{\beta}}f_{\alpha}f_{\beta
}\right] \\
=4bt^{3}\eta^{3}f_{0}\underset{\alpha,\beta=1}{\overset{n}{\sum}}\left[
\left(  A_{\alpha \beta}f_{\bar{\beta}\bar{\alpha}}+A_{\bar{\alpha}\bar{\beta}%
}f_{\beta \alpha}\right)  +\left(  A_{\alpha \beta,\bar{\alpha}}f_{\bar{\beta}%
}+A_{\bar{\alpha}\bar{\beta},\alpha}f_{\beta}\right)  +\left(  A_{\alpha \beta
}f_{\bar{\alpha}}f_{\bar{\beta}}+A_{\bar{\alpha}\bar{\beta}}f_{\alpha}%
f_{\beta}\right)  \right] \\
\geq-8bt^{3}\eta^{3}\left \vert f_{0}\right \vert \underset{\alpha,\beta
=1}{\overset{n}{\sum}}\left(  \left \vert A_{\bar{\alpha}\bar{\beta}%
}\right \vert \left \vert f_{\beta \alpha}\right \vert +\left \vert A_{\alpha
\beta,\bar{\alpha}}\right \vert \left \vert f_{\bar{\beta}}\right \vert
+\left \vert A_{\bar{\alpha}\bar{\beta}}\right \vert \left \vert f_{\alpha
}\right \vert \left \vert f_{\beta}\right \vert \right)
\end{array}
\label{16}%
\end{equation}
In (\ref{16}), by Young's inequality and noting that $t\leq1$, $\eta \leq1$, we
have following estimates:%
\begin{equation}%
\begin{array}
[c]{lll}%
-8bt^{3}\eta^{3}\left \vert f_{0}\right \vert \underset{\alpha,\beta=1}%
{\overset{n}{\sum}}\left \vert A_{\bar{\alpha}\bar{\beta}}\right \vert
\left \vert f_{\beta \alpha}\right \vert  & \geq & \underset{\alpha,\beta
=1}{\overset{n}{\sum}}-8k_{1}bt^{3}\eta^{3}\left \vert f_{0}\right \vert
\left \vert f_{\beta \alpha}\right \vert \\
& \geq & \underset{\alpha,\beta=1}{\overset{n}{\sum}}\left(  -4k_{1}bt^{3}%
\eta^{3}\left \vert f_{\beta \alpha}\right \vert ^{2}-4k_{1}bt^{3}\eta^{3}%
f_{0}^{2}\right) \\
& \geq & -4k_{1}bn^{2}\left(  t^{2}\eta^{2}f_{0}^{2}\right)  -4k_{1}bt^{2}%
\eta^{2}\underset{\alpha,\beta=1}{\overset{n}{\sum}}\left \vert f_{\beta \alpha
}\right \vert ^{2}%
\end{array}
\label{17}%
\end{equation}
and
\begin{equation}%
\begin{array}
[c]{lll}%
-8bt^{3}\eta^{3}\left \vert f_{0}\right \vert \underset{\alpha,\beta=1}%
{\overset{n}{\sum}}\left \vert A_{\alpha \beta,\bar{\alpha}}\right \vert
\left \vert f_{\bar{\beta}}\right \vert  & \geq & -8k_{1}bt^{3}\eta^{3}%
\underset{\alpha,\beta=1}{\overset{n}{\sum}}\left \vert f_{0}\right \vert
\left \vert f_{\bar{\beta}}\right \vert \\
& \geq & -8k_{1}bt^{3}\eta^{3}\underset{\alpha,\beta=1}{\overset{n}{\sum}%
}\left(  \frac{1}{2}f_{0}^{2}+\frac{1}{2}\left \vert f_{\bar{\beta}}\right \vert
^{2}\right) \\
& \geq & -4k_{1}bn^{2}t^{3}\eta^{3}f_{0}^{2}-4k_{1}bnt^{3}\eta^{3}%
\underset{\beta=1}{\overset{n}{\sum}}\left \vert f_{\bar{\beta}}\right \vert
^{2}\\
& \geq & -4k_{1}bn^{2}\left(  t^{2}\eta^{2}f_{0}^{2}\right)  -2k_{1}bn\left(
t\eta \left \vert \nabla_{b}f\right \vert ^{2}\right)
\end{array}
\label{18}%
\end{equation}
and%
\begin{equation}%
\begin{array}
[c]{lll}%
-8bt^{3}\eta^{3}\left \vert f_{0}\right \vert \underset{\alpha,\beta=1}%
{\overset{n}{\sum}}\left \vert A_{\bar{\alpha}\bar{\beta}}\right \vert
\left \vert f_{\alpha}\right \vert \left \vert f_{\beta}\right \vert  & \geq &
-8k_{1}bt^{3}\eta^{3}\left \vert f_{0}\right \vert \underset{\alpha,\beta
=1}{\overset{n}{\sum}}\left(  \frac{1}{2}\left \vert f_{\alpha}\right \vert
^{2}+\frac{1}{2}\left \vert f_{\beta}\right \vert ^{2}\right) \\
& \geq & -4k_{1}bt^{3}\eta^{3}\left \vert f_{0}\right \vert \left(
n\underset{\alpha=1}{\overset{n}{\sum}}\left \vert f_{\alpha}\right \vert
^{2}+n\underset{\beta=1}{\overset{n}{\sum}}\left \vert f_{\beta}\right \vert
^{2}\right) \\
& \geq & -4k_{1}bnt^{3}\eta^{3}\left \vert f_{0}\right \vert \left \vert
\nabla_{b}f\right \vert ^{2}\\
& \geq & -2k_{1}b^{2}nt^{3}\eta^{3}f_{0}^{2}\left \vert \nabla_{b}f\right \vert
^{2}-2k_{1}nt^{3}\eta^{3}\left \vert \nabla_{b}f\right \vert ^{2}\\
& = & -2k_{1}b^{2}n\left(  t\eta \left \vert \nabla_{b}f\right \vert ^{2}\right)
\left(  t^{2}\eta^{2}f_{0}^{2}\right)  -2k_{1}nt^{3}\eta^{3}\left \vert
\nabla_{b}f\right \vert ^{2}.
\end{array}
\label{19}%
\end{equation}

Substitute estimates (\ref{17}), (\ref{18}), and (\ref{19}) into (\ref{15}),
one obtains%
\[%
\begin{array}
[c]{lll}%
t\eta \Delta_{b}\left(  \eta F\right)  & \geq & \frac{1}{na^{2}}\left(  \eta
F\right)  ^{2}-\frac{C}{R}\left(  \eta F\right)  +2t\eta \left \langle
\nabla_{b}\eta,\text{ }\nabla_{b}F\right \rangle -2t\eta^{2}\left \langle
\nabla_{b}f,\text{ }\nabla_{b}F\right \rangle \\
&  & +4t^{2}\eta^{2}\left[  \left(  1-bk_{1}n\right)  \underset{\alpha
,\beta=1}{\overset{n}{\sum}}\left \vert f_{a\beta}\right \vert ^{2}%
+\underset{\alpha,\beta=1,\alpha \not =\beta}{\overset{n}{\sum}}\left \vert
f_{a\bar{\beta}}\right \vert ^{2}\right] \\
&  & +\left[  n-8bk_{1}n^{2}-\frac{bC}{R}-\left(  \frac{2b}{na^{2}}%
+2b^{2}k_{1}n+\frac{bC}{R}\right)  \left(  \eta F\right)  \right]  t^{2}%
\eta^{2}f_{0}^{2}\\
&  & +\left[  -\frac{2\left(  1+a\right)  }{na^{2}}\left(  \eta F\right)
-2k-2n\left(  1+b\right)  k_{1}-\frac{4}{b}\right]  t\eta \left \vert \nabla
_{b}f\right \vert ^{2}.
\end{array}
\]

Next as shown in the same computation as in Proposition \ref{prop3}, at the
maximal point $\left(  p(t),\text{ }t\right)  $%
\begin{equation}%
\begin{array}
[c]{lll}%
0 & \geq & \left(  \frac{1}{na^{2}}-\frac{C}{R}\right)  \left(  \eta F\right)
^{2}-\frac{3C}{R}\left(  \eta F\right) \\
&  & +4t^{2}\eta^{2}\left[  \left(  1-bk_{1}n\right)  \underset{\alpha
,\beta=1}{\overset{n}{\sum}}\left \vert f_{a\beta}\right \vert ^{2}%
+\underset{\alpha,\beta=1,\alpha \not =\beta}{\overset{n}{\sum}}\left \vert
f_{a\bar{\beta}}\right \vert ^{2}\right] \\
&  & +\left[  n-8bk_{1}n^{2}-\frac{bC}{R}-\left(  \frac{2b}{na^{2}}%
+2b^{2}k_{1}n+\frac{bC}{R}\right)  \left(  \eta F\right)  \right]  t^{2}%
\eta^{2}f_{0}^{2}\\
&  & +\left[  -\frac{2\left(  1+a\right)  }{na^{2}}\left(  \eta F\right)
-2k-2n\left(  1+b\right)  k_{1}-\frac{4}{b}-\frac{C}{R}\right]  t\eta
\left \vert \nabla_{b}f\right \vert ^{2}.
\end{array}
\label{20}%
\end{equation}

We claim at $t=1,$ there exists a small constant $b_{0}=b_{0}(n,k,k_{1})>0$
such that for any $0<b\leq b_{0}$
\[
\left(  \eta F\right)  \left(  p\left(  1\right)  ,\text{ }1,\text{ }R,\text{
}b\right)  <\frac{na^{2}}{-2\left(  1+a\right)  }\left(  2k+2n\left(
1+b\right)  k_{1}+\frac{4}{b}+\frac{C}{R}\right)
\]
if $R$ is large enough which to be determined later. Here $\left(  1+a\right)
<0$ for some $a$ to be chosen later (say $1+a=-\frac{5}{n}).$

We prove it by contradiction. Suppose not, that is
\[
\left(  \eta F\right)  \left(  p\left(  1\right)  ,\text{ }1,\text{ }R,\text{
}b\right)  \geq \frac{na^{2}}{-2\left(  1+a\right)  }\left(  2k+2n\left(
1+b\right)  k_{1}+\frac{4}{b}+\frac{C}{R}\right)  .
\]
Since $\left(  \eta F\right)  \left(  p\left(  t\right)  ,\text{ }t,\text{
}R,\text{ }b\right)  $ is continuous in the variable $t$ and $\left(  \eta
F\right)  \left(  p\left(  0\right)  ,\text{ }0,\text{ }R,\text{ }b\right)
=0$, by Intermediate-value theorem there exists a $t_{0}\in \left(  0,\text{
}1\right]  $ such that%
\[
\left(  \eta F\right)  \left(  p\left(  t_{0}\right)  ,\text{ }t_{0},\text{
}R,\text{ }b\right)  =\frac{na^{2}}{-2\left(  1+a\right)  }\left(
2k+2n\left(  1+b\right)  k_{1}+\frac{4}{b}+\frac{C}{R}\right)  .
\]

Now we apply (\ref{20}) at the point $\left(  p\left(  t_{0}\right)
,t_{0}\right)  $, denoted by $\left(  p_{0},t_{0}\right)  $. We have
\begin{equation}%
\begin{array}
[c]{l}%
\left[  \left(  \frac{1}{na^{2}}-\frac{C}{R}\right)  \left(  \eta F\right)
\left(  p_{0},t_{0}\right)  -\frac{3C}{R}\right] \\
=\left[  \left(  \frac{1}{na^{2}}-\frac{C}{R}\right)  \left(  \frac{na^{2}%
}{-2\left(  1+a\right)  }\right)  \left(  2k+2n\left(  1+b\right)  k_{1}%
+\frac{4}{b}+\frac{C}{R}\right)  -\frac{3C}{R}\right] \\
=\left \{  \frac{-1}{2(1+a)}\left(  2k+2n\left(  1+b\right)  k_{1}+\frac{4}%
{b}\right)  -\frac{C}{R}\left[  \frac{na^{2}}{-2\left(  1+a\right)  }\left(
2k+2n\left(  1+b\right)  k_{1}+\frac{4}{b}+\frac{C}{R}\right)  +\frac
{1}{2\left(  1+a\right)  }+3\right]  \right \}
\end{array}
\label{14}%
\end{equation}
and%

\begin{equation}%
\begin{array}
[c]{l}%
\left[  n-8bk_{1}n^{2}-\frac{bC}{R}-\left(  \frac{2b}{na^{2}}+2b^{2}%
k_{1}n+\frac{bC}{R}\right)  \left(  \eta F\right)  \left(  p_{0},t_{0}\right)
\right] \\
=n-8bk_{1}n^{2}-\frac{bC}{R}-\left(  \frac{2b}{na^{2}}+2b^{2}k_{1}n+\frac
{bC}{R}\right)  \left(  \frac{na^{2}}{-2\left(  1+a\right)  }\right)  \left(
2k+2n\left(  1+b\right)  k_{1}+\frac{4}{b}+\frac{C}{R}\right) \\
=n-8bk_{1}n^{2}-\frac{bC}{R}+(\frac{na^{2}}{2\left(  1+a\right)  })(\frac
{2b}{na^{2}}+2b^{2}k_{1}n)\left(  2k+2n\left(  1+b\right)  k_{1}+\frac{4}%
{b}+\frac{C}{R}\right) \\
+\frac{bC}{R}\left(  \frac{na^{2}}{2\left(  1+a\right)  }\right)  \left(
2k+2n\left(  1+b\right)  k_{1}+\frac{4}{b}+\frac{C}{R}\right) \\
=\{n-8bk_{1}n^{2}+(\frac{b+a^{2}b^{2}n^{2}k_{1}}{\left(  1+a\right)
})[2k+2n\left(  1+b\right)  k_{1}+\frac{4}{b}]\} \\
+\frac{C}{R}\{-b+(\frac{b+a^{2}b^{2}n^{2}k_{1}}{\left(  1+a\right)  }%
)+\frac{na^{2}b}{2(1+a)}[2k+2n\left(  1+b\right)  k_{1}+\frac{4}{b}+\frac
{C}{R}]\}.
\end{array}
\label{14a}%
\end{equation}

Now we choose $a$ and $b$ such that%
\begin{equation}%
\begin{array}
[c]{l}%
n-8bk_{1}n^{2}+(\frac{b+a^{2}b^{2}n^{2}k_{1}}{\left(  1+a\right)
})[2k+2n\left(  1+b\right)  k_{1}+\frac{4}{b}]\\
=n-b\{8k_{1}n^{2}-(\frac{1+a^{2}bn^{2}k_{1}}{\left(  1+a\right)
})[2k+2n\left(  1+b\right)  k_{1}]-(\frac{4a^{2}n^{2}k_{1}}{1+a})\}+\frac
{4}{1+a}\\
>0.
\end{array}
\label{14b}%
\end{equation}
This can be done by choosing
\[
(1+a)<-\frac{4}{n}%
\]
and then choose a small $b_{0}=b_{0}(n,k,k_{1})>0$ such that for any $b\leq
b_{0}$%
\[
n-b\{8k_{1}n^{2}-(\frac{1+a^{2}bn^{2}k_{1}}{\left(  1+a\right)  }%
)[2k+2n\left(  1+b\right)  k_{1}]-(\frac{4a^{2}n^{2}k_{1}}{1+a})\}+\frac
{4}{1+a}>0
\]
and
\[
\left(  1-bk_{1}n\right)  >0.
\]
In particular, we let
\[
1+a=-\frac{5}{n}.
\]
Then for any $0<b\leq b_{0}$, one obtains
\[
\left[  \left(  \frac{1}{na^{2}}-\frac{C}{R}\right)  \left(  \eta F\right)
\left(  p_{0},t_{0}\right)  -\frac{3C}{R}\right]  >0
\]
and
\[
\left[  n-8bk_{1}n^{2}-\frac{bC}{R}-\left(  \frac{2b}{na^{2}}+2b^{2}%
k_{1}n+\frac{bC}{R}\right)  \left(  \eta F\right)  \left(  p_{0},t_{0}\right)
\right]  >0
\]
for $R=R(b,k,k_{1})$ large enough.\ This leads to a contradiction with
(\ref{20}). Hence
\[
\left(  \eta F\right)  \left(  1,\text{ }p\left(  1\right)  ,\text{ }R,\text{
}b\right)  <\frac{na^{2}}{-\left(  1+a\right)  }\left(  k+n\left(  1+b\right)
k_{1}+\frac{2}{b}+\frac{C}{R}\right)  .
\]
This implies for $1+a=-\frac{5}{n}$
\[
\underset{x\in B\left(  2R\right)  }{\max}\left(  \left \vert \nabla
_{b}f\right \vert ^{2}+b\eta f_{0}^{2}\right)  \left(  x\right)  <\frac
{(n+5)^{2}}{5}\left(  k+n\left(  1+b\right)  k_{1}+\frac{2}{b}+\frac{C}%
{R}\right)  .
\]
When we fix on the set $x\in B\left(  R\right)  $, we obtain%
\[
\left \vert \nabla_{b}f\right \vert ^{2}+bf_{0}^{2}<\frac{(n+5)^{2}}{5}\left(
k+n\left(  1+b\right)  k_{1}+\frac{2}{b}+\frac{C}{R}\right)
\]
on $B\left(  R\right)  $. Note that the preceding computation is not valid if
$\eta F$ is not smooth at $x_{0}$. In this case, we may use a trick due to E.
Calabi ( see \cite{w} for details).

This completes the proof of Theorem \ref{t3}.
%TCIMACRO{\TeXButton{endproof}{\hfill}}%
%BeginExpansion
\hfill
%EndExpansion

\end{document}